\documentclass[reqno,11pt]{amsart}
\usepackage{amsthm}
\usepackage{amsmath,mathrsfs}
\usepackage{color}
\usepackage{amssymb}
\usepackage{hyperref}
\usepackage{enumerate}
\usepackage{latexsym,array}
\usepackage{amsfonts}
\usepackage{shadow}

\newtheorem{Pa}{Paper}[section]
\newtheorem{Tm}[Pa]{{\bf Theorem}}
\newtheorem{La}[Pa]{{\bf Lemma}}
\newtheorem{Dn}[Pa]{{\bf Definition}}
\newtheorem{Cy}[Pa]{{\bf Corollary}}

\newtheorem{Rk}[Pa]{{\bf Remark}}
\newtheorem{Pn}[Pa]{{\bf Proposition}}

\newtheorem{Ex}[Pa]{{\bf Example}}

\def\C{\mathbb C}
\def\hh{\mathbb{H}}
\author[D. Alpay]{Daniel Alpay}
\address{(DA) Department of Mathematics\\
Ben-Gurion University of the Negev\\
Beer-Sheva 84105 Israel} \email{dany@math.bgu.ac.il}

\author[F. Colombo]{Fabrizio Colombo}
\address{(FC) Politecnico di
Milano\\Dipartimento di Matematica\\Via E. Bonardi, 9\\20133
Milano, Italy}
\email{fabrizio.colombo@polimi.it}

\author[I. Lewkowicz]{Izchak Lewkowicz}
\address{(IL) Department of electrical engineering\\
Ben-Gurion University of the Negev\\
Beer-Sheva 84105 Israel} \email{izchak@ee.bgu.ac.il}

\author[I. Sabadini]{Irene Sabadini}
\address{(IS) Politecnico di
Milano\\Dipartimento di Matematica\\Via E. Bonardi, 9\\20133
Milano, Italy}
\email{irene.sabadini@polimi.it}
\title[]
{Realizations of slice hyperholomorphic generalized contractive
and positive functions} \oddsidemargin 0.2in \evensidemargin
0.2in \topmargin -0.5in \textwidth 15.8truecm \textheight 23.3truecm
\def\H{\mathbb H}

\def\R{\mathbb R}

\def\C{\mathbb C}

\def\(s){\mathscr S(\R\times\R)}



\thanks{D. Alpay thanks the Earl Katz family for endowing the chair
which supported his research, and the Binational Science
Foundation Grant number 2010117. F. Colombo and
I. Sabadini
acknowledge the Center for Advanced Studies of
the Mathematical Department of the Ben-Gurion University of the Negev for
the support and the kind hospitality during
the period in which part of this paper has been written.}
\begin{document}
\maketitle \tableofcontents
\parindent 0cm
\begin{abstract}
We introduce generalized Schur functions and generalized positive
functions in the setting of slice hyperholomorphic functions and
study their realizations in terms of associated reproducing kernel
Pontryagin spaces.  To this end, we also prove some results in quaternionic functional analysis like an invariant subspace theorem for contractions in a Pontryagin space.
We also consider slice hyperholomorphic functions on the half space  $\mathbb H_+$ of quaternions with positive real parts
and we study the Hardy space $\mathbf{H}_2(\mathbb H_+)$ and Blaschke products in this framework.
\end{abstract}

\parindent 0cm

\section{Introduction}
\setcounter{equation}{0}
In this paper we continue the study of Schur analysis in the
hyperholomorphic setting, initiated in \cite{acs1}, and continued
in \cite{acs3,acs2,abcs-np1}. To set the framework we
first recall a few facts on the classical case.

\subsection{Schur analysis}
Functions analytic and contractive in the open unit disk, or in an
open half-plane, play an important role in operator theory, signal
processing and related fields. Their study, and the study of their
counterparts in various settings, may be called Schur analysis;
see \cite{hspnw,MR2002b:47144,C-schur}. In the case of
matrix-valued, or operator-valued functions, contractivity is
considered with respect to an indefinite metric. An important
example is that of the characteristic operator function and
associated operator models.\smallskip

More precisely, let $T$ be a (say, bounded, for the present
discussion) self-adjoint operator in a Hilbert space such that
$T-T^*$ has finite rank, say $m$. Let
\[
T-T^*=\frac{CJC^*}{2i}
\]
where $J$ is a $m\times m$ matrix which is both self-adjoint and
unitary (a {\sl signature matrix}). Then, the {\sl matrix-valued}
function
\[
\Theta(z)=I+2iC^*(zI-T)^{-1}CJ
\]
is such that
\begin{equation}
\label{jcont}
\Theta(z)J\Theta(z)^*\ge J,
\end{equation}
for $z$ in the intersection $\Omega(T)$ of the upper open
half-plane and of the spectrum of $T$. The function $\Theta$ is
the characteristic operator function of $T$.\smallskip

Property \eqref{jcont} is called $J$-expansivity (or $-J$-contractivity), and is in
fact equivalent to the fact that the kernel
\begin{equation}
\label{ktheta}
K_\Theta(z,w)\stackrel{\rm def.}{=}
\frac{\Theta(z)J\Theta(w)^*-J}{-2i(z-w^*)}=C^*(zI-T)^{-1}(wI-T)^{-*}C
\end{equation}
is positive definite in $\Omega(T)$.
\smallskip

The function $\Theta$ provides a functional model for $T$, see
\cite{MR48:904}. A key fact in the theory is the multiplicative
structure of $J$-expansive functions, due to V. Potapov, see
\cite{pootapov}. We also refer to the historical note of M. Livsic \cite{MR1473265}.
It is also worth mentioning the original papers of M. Livsic \cite{MR0049488,livsic-54},
where the notion of
characteristic operator
function first appears.\smallskip

Replacing $J$ by $-J$ we obtain $J$-contractive, rather than
$J$-expansive functions, and this is the choice we make in the
sequel (see in particular the last section).

\subsection{Negative squares}
One can consider functions $\Theta$ such that the associated
kernel $K_\Theta$ has a finite number of negative squares (see
Definition \ref{def:neg}), rather than being positive definite.
Such classes of operator-valued functions were introduced and
studied by Krein and Langer in a long series of papers, see for
instance \cite{kl2,MR518342,MR563344,MR614775,kl3}.
These works are set in the framework of the open upper
half-plane. To make a better connection with the quaternionic
case we consider here the open right half-plane. In the complex
variable case, the two cases are equivalent via a conformal map.
This differs from the quaternionic setting, as will be clear
in the sequel. First recall that an operator $J$ in a Hilbert
space $\mathcal H$ is called a signature operator if it is
self-adjoint and unitary. Its spectrum is then concentrated on
$\pm 1$. When $-1$ is an eigenvalue of finite order, we denote
this multiplicity by $\nu_-(J)$. Let $\mathcal H_1$ and $\mathcal
H_2$ be two Hilbert spaces and let $J_1\in\mathbf L(\mathcal
H_1)$ and $J_2\in\mathbf L(\mathcal H_2)$ be two signature
operators, such that $\nu_-(J_1)=\nu_-(J_2)<\infty$. The $\mathbf
L(\mathcal H_1,\mathcal H_2)$-valued function $S$ analytic in an
open subset $\Omega$ of the  open right half-plane is called a
generalized Schur function if the kernel
\begin{equation}
\label{schurgeneral}
\frac{J_2-S(z)J_1S(w)^*}{z+\overline{w}}
\end{equation}
has a finite number of negative squares in $\Omega$.\\

For instance when
$J_1=J_2=\begin{pmatrix}0&-i\\i&0\end{pmatrix}$, functions in the
corresponding class are introduced and used in \cite{MR518342} to
describe the set of all generalized resolvents of a given
Hermitian operator, see \cite[Satz 3.5, p. 407 and Satz 3.9, p.
409]{MR518342}.\smallskip

Similarly given a Hilbert space $\mathcal H$ and a signature
operator $J$ (possibly with $\nu_-(J)=\infty$, see \cite[p. 358,
footnote]{kl1}), a $\mathbf L(\mathcal H)$-valued function $\Phi$
analytic in some open subset $\Omega$ of the right open
half-plane $\Pi_+$ is called generalized positive if the kernel
\begin{equation}
\label{kphi}
\frac{J\Phi(z)+\Phi(w)^*J}{z+w^*}
\end{equation}
has a finite number of negative squares in $\Omega$. The
$Q$-function of an Hermitian operator in a Pontryagin space,
introduced by Krein and Langer has such a property, see
\cite{MR47:7504}. The function $\Phi$ will be called positive if
the kernel \eqref{kphi} is positive definite.
\smallskip

In both cases, Krein and Langer proved in the above mentioned
works, among numerous results, realization formulas which ensure
the existence of a meromorphic extension to the whole of $\Pi_+$.
It is worth mentioning that a key result to prove this extension
is that the part of the spectrum of a contraction in a Pontryagin
space which lies outside the closed unit disk consists only of a
finite number of eigenvalues. One proof of this fact uses the
Schauder-Tychonoff fixed point theorem (see the discussion
\cite[p. 248]{MR92m:47068}). We also mention that a study of
generalized Schur function of the open unit disk has been given
in \cite{adrs} and that unified formulas for a number of cases
which include the line and circle case were developed in
\cite{abds2-jfa}, based on an approach including both the disk
and half-plane cases developed in \cite{abds2-jfa,ad3}.\\

Finally we mention the works \cite{MR2828331,MR2454941,MR1821917}
to stress the interest of positive and generalized positive
functions in linear system theory and operator theory.

\subsection{The slice hyperholomorphic case}

In previous papers we extended results of Schur analysis in the
slice hyperholomorphic case, in the setting of the unit ball
$\mathbb B_1$ of the quaternions. We considered in \cite{acs1}
the Schur algorithm, and the underlying counterpart of the Hardy
space. Blaschke products and related interpolation problems in the
Hardy space were studied in \cite{acs2}. Nevanlinna-Pick
interpolation for Schur functions is studied in \cite{abcs-np1}
while the case of kernels having a number of negative squares was
studied in \cite{acs3}.\smallskip

In contrast to the above mentioned papers, we consider in this
work functions which are slice hyperholomorphic in an open subset
of the open half-space
\[
\mathbb H_+=\left\{p\in\mathbb H\,\,;{\rm Re}\, p>0\right\},
\]
which intersects the positive real axis.\smallskip

We define and study the counterparts of the kernel
\eqref{schurgeneral} and \eqref{kphi} in the setting of slice
hyperholomorphic functions. Here we consider the case of
operator-valued generalized positive functions and generalized
Schur functions, rather than scalar or matrix-valued functions.
The extension of realization of generalized positive
functions to the slice hyperholomorphic setting,
introduced in this work, calls upon a corresponding
extension of the Kalman-Yakubovich-Popov lemma (also known as  Positive Real Lemma;  see the discussion for the classical case in the next paragraph). This will be addressed in another work.\smallskip

A $\mathbb C^{p\times p}$-valued
function $~F(s)$, analytic in $\mathbb C_+$ is said to be
{\em positive} if
\begin{equation}\label{eq:PosFunc}
\begin{matrix}
F(s)+F(s)^*\ge 0 &~&s\in\C_+,\end{matrix}
\end{equation}
where the inequality sign means that the Hermitian matrix is non negative, and where
$\C_+$ denotes the open right half of the
complex plane. The study of rational positive functions has been motivated from the 1920's by (lumped)
electrical networks theory, see e.g. \cite{AV}, \cite{Be}. From
the 1960's positive functions also appeared in books on absolute
stability theory, see e.g.
\cite{Po}. A $\mathbb C^{p\times p}$-valued function of bounded type in $\mathbb C_+$
(i.e. a quotient of two functions analytic and bounded in
$\mathbb C_+$) is called  generalized positive if
\begin{equation}\label{DefGP}
F(i\omega)+F(i\omega)^*\ge 0,
\quad a.e.\quad\quad\omega\in\R,
\end{equation}
where $F(i\omega)$ denotes the non-tangential
limit\begin{footnote}
{This limit exists almost everywhere on $i\R$ because $F$ is
assumed of bounded type in $\mathbb C_+$, see e.g.
\cite{Du}.}
\end{footnote}
of $F$ at the point $i\omega$.
In the classical setting, generalized positive functions were introduced in the context
of the Positive Real Lemma (PRL), see \cite{AM} and references
therein\begin{footnote}{The original formulation was real.
The case we address is in fact {\em generalized} positive and
{\em complex}, but we wish to adhere to the commonly used
term: Positive Real Lemma.}\end{footnote}. For applications of
generalized positive functions
see \cite{HSK}.
The renowned Kalman-Yakubovich-Popov Lemma, which has been recognized as a
fundamental result in System Theory,
establishes a connection between two presentations
of positive functions, as rational functions
and the respective state space realization, see e.g. \cite{AV}, \cite{MR525380}.
For its extension to generalized positive functions,
see \cite{AM}, \cite{DDGK}.\\

The paper consists of seven sections, besides the Introduction. In
Section 2 we recall the notion of quaternionic Pontryagin
spaces,  and we discuss some preliminaries on negative squares,
kernels and realizations; then we provide some preliminaries on
slice hyperholomorphic functions, the class of functions that we
use in this paper. Section 3 deals with operator-valued slice hyperholomorphic functions, their products and a
useful property of extension, see Proposition 3.24. In Section 4 we study the Hardy space of the
half-space $\mathbb H_+$ of quaternions with real
positive part and Blaschke factors and products in this framework. Then, in Section 5 we provide the proof of the
quaternionic version of the Schauder-Tychonoff fixed point
theorem whose proof is not substantially different from the one
in the complex case, but we insert it for the sake of
completeness. This result is crucial to show an invariant
subspace theorem for contractions in Pontryagin spaces. Sections 6
and 7 deal with the study of kernels with a finite number of negative squares and associated with generalized Schur functions and we prove a realization theorem in this setting.  We also give as an example the characteristic operator
function of a quaternionic non anti-self-adjoint operator.
Section 8 deals with realizations for generalized positive
functions. We also define the positive function associated
to a pair of anti-self-adjoint operators. The properties of this
function will be presented in a future publication.
\section{Preliminaries}
\setcounter{equation}{0}

In this section, which is divided into three subsections,  we
collect a number of facts respectively  on Pontryagin spaces,
slice hyperholomorphic functions and their realizations.
\subsection{Negative squares and kernels}
\setcounter{equation}{0}

An important role in this paper is played by quaternionic Pontryagin spaces, and we
first recall this notion. We refer to \cite{acs4,as3} for more details. Let $\mathcal V$ be
a right quaternionic vector space endowed with a Hermitian form (also called
inner product) $[\cdot,\cdot]$ from $\mathcal V\times \mathcal V$ into $\mathbb H$,
meaning that:
\[
\begin{split}
[ua+vb,w]&=[u,w]a+[v,w]b,\\
[v,w]&=\overline{[w,v]},
\end{split}
\]
for all choices of $u,v,w\in\mathcal V$ and $a,b\in\mathbb H$.
In particular the inner product $[\cdot,\cdot]$ satisfy
\[
[va, wb]=\overline{b}[ v,w] a.
\]
When the space $\mathcal V$ is two-sided, we require that
\begin{equation}
[ v, aw]=[ \overline{a}v, w],\quad a\in\mathbb H,\quad
v,w\in\mathcal V.
\label{rightcondition}
\end{equation}
Condition \eqref{rightcondition} is used in particular in the proof of formula
\eqref{I3}.\\

\begin{Dn}
The space $\mathcal V$ is called a right-quaternionic Pontryagin space if there exists two subspaces
$\mathcal V_+$, $\mathcal V_-$ of $\mathcal V$ such that $\mathcal V=\mathcal V_++\mathcal V_-$ and:\\
$(i)$ The space $\mathcal V_+$ endowed with $[\cdot,\cdot]$ is
a  right-quaternionic Hilbert space.\\
$(ii)$ The space $\mathcal V_-$ endowed with $-[\cdot,\cdot]$ is
a finite dimensional right-quaternionic Hilbert space.\\
$(iii)$ The sum $\mathcal V_++\mathcal V_-$ is direct and orthogonal, meaning that $\mathcal V_+\cap\mathcal V_-=
\left\{0\right\}$ and $[v_+,v_-]=0$ for every choice of $v_+\in\mathcal V_+$ and
$v_-\in\mathcal V_-$.
\end{Dn}

We denote a direct and orthogonal sum by
\begin{equation}
\label{decomportho}
\mathcal V=\mathcal V_+[\stackrel{\cdot}{+}]\mathcal V_-.
\end{equation}
In general, such a decomposition will not be  unique. The
inner product
\[
\langle v,w\rangle=[v_+,w_+]-[v_-,w_-]
\]
where $v_\pm, w_\pm\in\mathcal V_\pm$,
makes $\mathcal V$ into a Hilbert space. The inner product
depends on the decomposition, but all the associated topologies
are equivalent. We refer to  \cite{as3} for more details on these
facts in the quaternionic case, while the case of the field of
complex numbers we refer to \cite{bognar}.\\

We now recall a few facts on matrices with quaternionic entries and on kernels,
which we will need in the sequel. A matrix $A\in\mathbb H^{m\times m}$ can be
written in a unique way as
\[
A=A_1+A_2j,
\]
where $A_1$ and $A_2$ belong to $\mathbb C^{m\times m}$. The map $\chi :\, \mathbb C^{m\times m} \to \mathbb C^{2m\times 2m}$ defined by
\begin{equation}
\label{defchi}
\chi(A)=\begin{pmatrix}A_1&A_2\\-\overline{A_2}&\overline{A_1}\end{pmatrix}
\end{equation}
satisfies
\[
\chi(AB)=\chi(A)\chi(B)\quad{\rm and}\quad \chi(A^*)=(\chi(A))^*.
\]
See for instance  \cite[Theorem 4.2, p. 29]{MR97h:15020}
(see also \cite[Proposition 3.8, p. 439]{as3}). The result itself is due to
Lee \cite{MR12:153i}.\smallskip

A key fact is that $A\in\mathbb H^{m\times m}$ is Hermitian (that is,
$A=A^*$) if and only if it can be written as $UDU^*$, where
$U\in\mathbb H^{m\times m}$ is unitary and $D\in\mathbb R^{m\times m}$ is
diagonal. The matrix $D$ is uniquely determined up to permutations, and one can define the signature of
an Hermitian matrix with quaternionic entries as the signature of $D$, see
\cite[Corollary 6.2, p. 41]{MR97h:15020} and the references
therein. The following result follows from the properties of $\chi$
and can be found in \cite[Proposition 3.16, p. 442]{as3}.

\begin{La}
\label{la:signa}
Assume $A\in\mathbb H^{m\times m}$ Hermitian. Then $A$ has signature
$(\nu_+,\nu_-,\nu_0)$ if and only if $\chi(A)$ has signature
$(2\nu_+,2\nu_-,2\nu_0)$.
\end{La}

We now turn to the notion of kernels having a finite number of negative squares.

\begin{Dn}
\label{def:neg}
Let $\mathcal H$ be a two-sided quaternionic
Hilbert space, with inner product $\langle\cdot,\cdot\rangle$,
and let $K(z,w)$ be a $\mathbf L(\mathcal H,\mathcal H)$-valued
function defined for $z,w$ in some set $\Omega$. The kernel is
called Hermitian if
\[
K(z,w)=K(w,z)^*,\quad z,w\in\Omega.
\]
It is said to have a finite number (say $\kappa$) of negative squares if
for every choice of $N\in\mathbb N$, of vectors
$c_1,\ldots, c_N\in \mathcal H$ and of points $w_1,\ldots, w_N\in\Omega$,
the $N\times N$ Hermitian matrix with $(u,v)$ entry
\[
[K(w_u,w_v)c_u,c_v]
\]
has at most $\kappa$ strictly negative eigenvalues, counted with multiplicities, and exactly
$\kappa$ strictly negative eigenvalues for some choice of $N,
w_1,\ldots, w_N$ and $c_1,\ldots, c_N$.
\end{Dn}

When $\kappa=0$ we have the
classical notion of positive definite function. Given a set $\Omega$, the one-to-one
correspondence between positive definite functions on $\Omega$ and reproducing kernel
Hilbert spaces of functions defined on $\Omega$ extends to a one-to-one correspondence
between functions having a finite number of negative squares and reproducing kernel
Pontryagin spaces (for more information on these spaces see \cite{as3}). This fact is due to P. Sorjonen \cite{Sorjonen73}
and L. Schwartz \cite{schwartz} in the complex case,
and is proved in \cite{as3} in the quaternionic case.\\
We conclude by mentioning a result,  \cite[Proposition 5.3]{acs3}, which will be used in the sequel:
\begin{Pn}
\label{pn51} Assume that $K(p,q)$  is $\mathbb H^{N\times
N}$-valued and has $\kappa$ negative squares in $V$ and let
$\alpha(p)$ be a $\mathbb H^{N\times N}$-valued slice
hyperholomorphic function and such that $\alpha(0)$ is
invertible. Then the function
\begin{equation}
\label{aka}
B(p,q)=\alpha(p)\star K(p,q)\star_r\alpha(q)^*
\end{equation}
has $\kappa$ negative squares in $V$.
\end{Pn}

\subsection{Slice hyperholomorphic functions}
Let $\hh$ be the real associative algebra of quaternions
with respect to the basis $\{1, i,j,k \}$
satisfying the relations
$
i^2=j^2=k^2=-1,\
 ij =-ji =k,\
jk =-kj =i ,
 \  ki =-ik =j .
$
We will denote a quaternion $p$ as $p=x_0+ix_1+jx_2+kx_3$,
$x_i\in \mathbb{R}$, its conjugate as
$\bar p=x_0-ix_1-jx_2-kx_3$, and $|p|^2=p\overline{p}$.
The real part $x_0$ of a quaternion will be denoted also by ${\rm Re}(p)$,
 $\mathbb{S}$ is the 2-sphere of purely imaginary unit quaternions, i.e.
$$
\mathbb{S}=\{ p=ix_1+jx_2+kx_3\ |\  x_1^2+x_2^2+x_3^2=1\}.
$$
Note that if $I\in\mathbb S$ then $I^2=-1$ and  a nonreal quaternion $p=x_0+ix_1+jx_2+kx_3$ uniquely determines an element $I_p=ix_1+jx_2+kx_3/|ix_1+jx_2+kx_3|\in\mathbb S$. When $p$ is real, then $p=p+I0$ for all $I\in\mathbb S$.
\begin{Dn}
Given $p\in\hh$, $p=p_0+I_pp_1$ we denote by $[p]$ the set of all elements of the form $p_0+Jp_1$ when $J$ varies in $\mathbb{S}$.
\end{Dn}
\noindent
The set $[p]$ is a $2$-sphere (we will often write that $[p]$ is  a sphere, for short) which is reduced to the point $p$ when $p\in\mathbb{R}$.
We now recall the definition of slice hyperholomorphic function.
\begin{Dn}
Let $\Omega\subseteq\hh$ be an open set and let
$f:\ \Omega\to\hh$ be a real differentiable function. Let
$I\in\mathbb{S}$ and let $f_I$ be the restriction of $f$ to the
complex plane $\mathbb{C}_I := \mathbb{R}+I\mathbb{R}$ passing through $1$
and $I$ and denote by $x+Iy$ an element on $\mathbb{C}_I$.
 We say that $f$ is a left slice hyperholomorphic (or slice hyperholomorphic or slice regular) function
  in $\Omega$ if, for every
$I\in\mathbb{S}$, we have
$$
\frac{1}{2}\left(\frac{\partial }{\partial x}+I\frac{\partial
}{\partial y}\right)f_I(x+Iy)=0.
$$
 We say that $f$ is a right slice hyperholomorphic function
 in $\Omega$ if, for every
$I\in\mathbb{S}$, we have
$$
\frac{1}{2}\left(\frac{\partial }{\partial x}f_I(x+Iy)+\frac{\partial
}{\partial y}f_I(x+Iy) I\right)=0.
$$
\end{Dn}

Slice hyperholomorphic functions have a nice behavior on the so called axially symmetric  slice domains defined below.
\begin{Dn}
Let $\Omega$ be a domain in $\mathbb{H}$.
We say that $\Omega$ is a
slice domain (s-domain for short) if $\Omega \cap \mathbb{R}$ is non empty and if
$\Omega\cap \mathbb{C}_I$ is a domain in $\mathbb{C}_I$ for all $I \in \mathbb{S}$.
We say that $\Omega$ is
axially symmetric if, for all $q \in \Omega$, the
sphere $[q]$ is contained in $\Omega$.
\end{Dn}
\begin{Rk}{\rm  Assume that $f:\ \Omega\subseteq \mathbb{C}\cong \mathbb{C}_I\to \mathbb H$ is a holomorphic map. Let $U_\Omega$ be the axially symmetric completion of $\Omega$, i.e. $U_\Omega=\bigcup_{J\in\mathbb S, x+Iy\in\Omega} (x+Jy)$. Its left slice hyperholomorphic extension ${\rm ext}(f): \ U_\Omega\subseteq \mathbb H \to\mathbb H$ is computed as follows (see \cite{MR2752913}):
\begin{equation}\label{ext}
{\rm ext}(f)(x+J y)=\frac 12 \left[ f(x+Iy) +f(x-Iy) + J I
(f(x-Iy)-f(x+Iy))\right].
\end{equation}
It is immediate that  ${\rm ext}(f+g)={\rm ext}(f)+{\rm ext}(g)$
and that if $f(z)=\sum_{n=0}^\infty f_n(z)$ then ${\rm
ext}(f)(z)=\sum_{n=0}^\infty {\rm ext}(f_n)(z)$.
It is also useful to recall that any  function $h$ slice hyperholomorphic on an axially symmetric s-domain $\Omega$ satisfies the formula, see \cite[Theorem 4.3.2]{MR2752913}
\begin{equation}\label{repr}
h(x+J y)=\frac 12 \left[ h(x+Iy) +h(x-Iy) + J I
(h(x-Iy)-h(x+Iy))\right].
\end{equation}
}
\end{Rk}

Let $f,g:\ \Omega \subseteq\mathbb{H}$ be slice hyperholomorphic functions.
Their restrictions to the complex plane $\mathbb{C}_I$ can be decomposed as
$f_I(z)=F(z)+G(z)J$,
$g_I(z)=H(z)+L(z)J$ where $J\in\mathbb{S}$, $J\perp I$ where $F$,
$G$, $H$, $L$ are holomorphic functions of the variable $z\in
\Omega \cap \mathbb{C}_I$, see \cite{MR2752913}, p. 117.
The $\star_l$-product of  $f$ and $g$, see \cite{MR2752913}, p. 125, is defined as the unique
left slice hyperholomorphic function whose restriction to the
complex plane $\mathbb{C}_I$ is given by
\begin{equation}\label{starproduct}
\begin{split}
(f_I\star_r g_I)(z):&=
(F(z)+G(z)J)\star_l(H(z)+L(z)J)\\
&=
(F(z)H(z)-G(z)\overline{L(\bar z)})+(G(z)\overline{H(\bar z)}+F(z)L(z))J.
\end{split}
\end{equation}
If $f,g$ are right slice hyperholomorphic, then with the above notations we have
$f_I(z)=F(z)+JG(z)$,
$g_I(z)=H(z)+JL(z)$
and
\begin{equation}\label{starproductright}
\begin{split}
(f_I\star_r g_I)(z):&=(F(z)+JG(z))\star_r(H(z)+JL(z))\\
&=
(F(z)H(z)-\overline{G(\bar z)}L(z))+J(G(z)H(z)+\overline{F(\bar z)}L(z))J,
\end{split}
\end{equation}
and $f\star_r g={\rm ext}(f_I\star_r g_I)$.

\begin{Rk}{\rm
In the sequel, we will consider functions $k(p,q)$ left slice hyperholomorphic in $p$ and right slice hyperholomorphic in $\bar q$. When taking the $\star$-product of a  function $f(p)$ slice hyperholomorphic in the variable $p$ with such a function $k(p,q)$, we will write $f(p)\star k(p,q)$ meaning that the $\star$-product is taken with respect to the variable $p$;  similarly, the $\star_r$-product of $k(p,q)$ with  functions right slice hyperholomorphic in the variable $\bar q$ is always taken with respect to $\bar q$.
}
\end{Rk}
Let $\Omega$ be an axially symmetric s-domain and let $p_0\in\Omega$.
Let us consider a function $f$ slice hyperholomorphic in $\Omega$ and assume that, in a neighborhood of $p_0$ in $\Omega$, it can be written in the form $f(p)=
\sum_{n=-\infty}^{+\infty} (p-p_0)^{\star n} a_n$ where $a_n\in
\mathbb H$.
\\

Following the standard nomenclature and
\cite{MR2955794} we have:
\begin{Dn}
A function $f$
  has a pole at the point $p_0$ if there exists $m\geq 0$ such that $a_{-k}=0$ for $k>m$.   The minimum of such $m$ is called the order of the pole;\\
  If $p$ is not a pole then we call it an essential singularity for $f$;\\
  $f$ has a removable singularity at $p_0$ if it can be  extended in a neighborhood of $p_0$ as a slice hyperholomorphic function.\\
\end{Dn}
     Note the following important fact: a function $f$ has a pole at $p_0$ if and only if its restriction to a complex plane has a pole. Note that there can be poles of order $0$: let us consider for example the function $(p+I)^{-\star}=(p^2+1)^{-1}(p-I)$. It has a pole of order $0$ at the point $-I$ which, however, is not a removable singularity, see \cite[p.55]{MR2752913} also for the definition of the $\star$-inverse.

\begin{Dn} Let $\Omega$ be an axially symmetric s-domain in $\mathbb H$. We say that a function $f:\, \Omega\to \mathbb H$ is
slice hypermeromorphic in $\Omega$ if $f$ is slice hyperholomorphic in
$\Omega'\subset \Omega$ such that $\Omega\setminus\Omega'$  has no point limit in $\Omega$ and every point in $\Omega\setminus\Omega'$ is a pole.
\end{Dn}
The functions which are slice hypermeromorphic are called
semi-regular in \cite{MR2955794} and for these functions we have
the following result, proved in \cite[Proposition 7.1, Theorem
7.3]{MR2955794}:
\begin{Pn}
Let $\Omega$ be an axially symmetric s-domain in $\mathbb H$
and let $f,g:\, \Omega\to\mathbb H$ be slice hyperholomorphic.
Then the function $f^{-\star}\star g$ is slice hypermeromorphic
in $\Omega$. Conversely, any slice hypermeromorphic function on
$\Omega$ can be locally expressed as $f^{-\star}\star g$ for
suitable $f$ and $g$.
\end{Pn}
\begin{Rk}{\rm
Since $f^{-\star}=(f\star f^c)^{-1}f^c$ (see \cite{MR2752913} for the notation)
it is then clear that the poles of a slice hypermeromorphic function occur in
correspondence to the zeros of the function $f\star f^c$ and so they are isolated
spheres, possibly reduced to real points.}
\end{Rk}

\section{Slice hyperholomorphic operator-valued functions}
\setcounter{equation}{0}
By begin the section by characterizing slice hyperholomorphic functions as those functions which
admit left derivative on each complex plane $\C_I$:
\begin{Dn}\label{slicederivative}
Let $f: \Omega\subseteq\H\to\H$ and let $p_0\in U$ be a nonreal point,
$p_0=u_0+Iv_0$. Let $f_I$ be the restriction of $f$ to the plane
$\C_I$. Assume that
\begin{equation}\label{deriv}
\lim_{p\to p_0, \,  p\in\C_I} (p-p_0)^{-1}(f_I(p)-f_I(p_0))
\end{equation}
exists. Then we say that $f$ admits left slice derivative in
$p_0$. If $p_0$ is real, assume that
\begin{equation}\label{deriv1}
\lim_{p\to p_0, \, p\in\C_I} (p-p_0)^{-1}(f_I(p)-f_I(p_0))
\end{equation}
exists, equal to the same value, for all $I\in\mathbb{S}$. Then we say that $f$ admits left
slice derivative in $p_0$. If $f$ admits left slice derivative
for every $p_0\in \Omega$ then we say that $f$ admits left slice
derivative in $\Omega$ or, for short, that $f$ is {\em left slice
differentiable} in $\Omega$.
\end{Dn}
It is possible to give an analogous definition for right slice differentiable functions: it is sufficient to multiply $(p-p_0)^{-1}$ on the right. In this case we will speak of right slice hyperhomolomorphic functions. In this paper,  we will speak of  slice differentiable functions or slice hyperholomorphic functions when we are considering them on the left, while we will specify if we consider the analogous notions on the right.\\
We have the following result:
\begin{Pn}\label{equivalence}
Let $\Omega\subseteq\H$ be an open set  and let $f: \Omega\subseteq\H\to\H$
be a real differentiable function. Then $f$ is slice
hyperholomorphic on $\Omega$ if and only if it admits slice derivative
on $\Omega$.
\end{Pn}
\begin{proof}
Let $f$ be a slice hyperholomorphic function on $\Omega$. Then its
restriction  to the complex plane $\C_I$  can be written as
$f_I(p)=F(p)+G(p)J$ where $J$ is any element in $\mathbb{S}$
orthogonal to $I$, $p$ belongs to $\C_I$ and $F,G:
\Omega\cap\C_I\to\C_I$ are holomorphic functions. Let $p_0$ be a
nonreal quaternion and let $p_0\in \Omega\cap\C_I$. Then we
have
\begin{equation}\label{limit}
\lim_{p\to p_0, \, p\in\C_I}
(p-p_0)^{-1}(f_I(p)-f_I(p_0))=\lim_{p\to p_0,\,  p\in\C_I}
(p-p_0)^{-1}(F(p)+G(p)J-F(p_0)-G(p_0)J)
\end{equation}
$$
=F'(p_0)+G'(p_0)J
$$
so the limit exists and $f$ admits slice derivative at every
nonreal point in $\Omega$. If $p_0$ is real then the same reasoning
shows that the limit in (\ref{limit}) exists on each complex
plane $\mathbb{C}_I$. Moreover, since $f$ is slice hyperholomorphic at $p_0$ we have
\[
F'(p_0)+G'(p_0)J=\frac 12\left(\frac{\partial}{\partial x}- I\frac{\partial}{\partial y}\right)(F+GJ) (p_0)= \frac{\partial}{\partial x} f(p_0)
\]
and so the limit exists on $\C_I$ for all $I\in\mathbb S$ equal to $\frac{\partial}{\partial x} f(p_0)$.
\\
Conversely, assume that $f$ admits slice derivative in $\Omega$. By
(\ref{deriv}) and (\ref{deriv1}) $f_I$ admits derivative on
$\Omega\cap \C_I$ for all $I\in\mathbb{S}$. Decomposing $f$ into complex components as
$f_I(p)=F(p)+G(p)J$, where $F,G: \Omega\cap\C_I\to\C_I$, $p=x+Iy$ and
$J$ is orthogonal to $I$, we deduce that both $F$ and $G$ admits complex derivative and thus they are
in the kernel of the Cauchy Riemann
operator $\partial_x+I\partial_y$ for all $I\in\mathbb{S}$ as well as $f_I$. Thus
$f$ is slice hyperholomorphic.
\end{proof}
\begin{Rk}{\rm
The terminology of Definition \ref{slicederivative} is consistent with the notion of slice derivative $\partial_s f$ of $f$, see \cite{MR2752913},
which is defined by:
\begin{displaymath}
\partial_s(f)(p) = \left\{ \begin{array}{ll}
\frac{1}{2}\left (\frac{\partial}{\partial x}f_{I}(x+ I y) - I\frac{\partial}{\partial y}f_{I}(x+ I y)\right ) & \textrm{ if $p=x+Iy$, \ $y\neq 0$},\\ \\
\displaystyle\frac{\partial f}{\partial x} (p) & \textrm{ if\
$p=x\in\mathbb{R}$.}
\end{array} \right.
\end{displaymath}
It is immediate that,  analogously to what happens in the complex
case, for any slice hyperholomorphic function we have $
\partial_s(f)(x+Iy) =\partial_x(f)(x+Iy)$. }
\end{Rk}
In the sequel, let $\mathcal{X}$ denote a left  quaternionic Banach space
and let $\mathcal{X}^*$ denote its dual, i.e. the set of
bounded, left linear maps from $\mathcal{X}$ to $\H$. In order to have that $\mathcal X^*$ has a structure of quaternionic
linear space, it is necessary to require that $\mathcal{X}$ is two sided quaternionic vector space. In this case, $\mathcal X^*$ turns out to be
a right vector space over $\mathbb H$.
\begin{Dn}
Let $\mathcal{X}$ be a two sided quaternionic Banach space and let $\mathcal{X}^*$ be its dual.
Let $\Omega$ be an open set in $\H$.\\
A function $f:\Omega\to \mathcal{X}$ is said to be {\em weakly slice hyperholomorphic}
in $\Omega$ if $\Lambda f$ admits slice derivative for every $\Lambda\in \mathcal{X}^*$.\\
A function $f:\Omega\to \mathcal{X}$ is said to be {\em strongly slice
hyperholomorphic} in $\Omega$ if
\begin{equation}\label{derivstrong}
\lim_{p\to p_0, p\in\C_I} (p-p_0)^{-1}(f_I(p)-f_I(p_0))
\end{equation}
exists in the topology of $\mathcal{X}$ in case $p_0\in \Omega$ is nonreal and
$p_0\in\C_I$ and if
\begin{equation}\label{deriv1strong}
\lim_{p\to p_0, p\in\C_I} (p-p_0)^{-1}(f_I(p)-f_I(p_0))
\end{equation}
exists in the topology of $\mathcal{X}$ for every $I\in\mathbb{S}$, equal to the same value, in case
$p_0\in \Omega$ is real.
\end{Dn}
Since the functionals $\Lambda\in \mathcal{X}^*$ are continuous, every
strongly slice hyperholomorphic function is weakly slice
hyperholomorphic. As it happens in the complex case, let us show
that also the converse is true.

To this end, let us observe that the following lemma holds. We
omit the proof since it works exactly as in the complex case (see
e.g. \cite{MR751959}, p. 189).
\begin{La}\label{techlemma}
Let $\mathcal{X}$ be a two sided quaternionic Banach space. Then a sequence
$\{v_n\}$ is Cauchy if and only if $\{\Lambda v_n\}$ is Cauchy
uniformly for $\Lambda\in \mathcal{X}^*$, $\|\Lambda\|\leq 1$.
\end{La}
\begin{Tm}
Every weakly slice hyperholomorphic function on $\Omega\subseteq\H$ is
strongly slice hyperholomorphic on $\Omega$.
\end{Tm}
\begin{proof}
The proof will follow the lines of the proof in the
complex case in \cite{MR751959}, p. 189. Let $f$ be a weakly
slice hyperholomorphic function on $\Omega$. Then, for any $\Lambda\in
\mathcal{X}^*$ and any $I\in\mathbb{S}$, we can choose $J\in\C_I$ and write
$(\Lambda f)_I(p)=(\Lambda f)_I(x+Iy)=F_\Lambda (x+Iy)+G_\Lambda
(x+Iy) J$ where $F_\Lambda, G_\Lambda : \C_I\to\C_I$. By
hypothesis, for any $p_0\in \Omega\cap \C_I$ the limit $\lim_{p\to
p_0,\, p\in\C_I} (p-p_0)^{-1}((\Lambda f)_I(p)- (\Lambda f)_I (p_0))$ exists, and so the
limits
$$
\lim_{p\to p_0,\, p\in\C_I}
(p-p_0)^{-1}(F_\Lambda(p)-F_{\Lambda}(p_0))\qquad \lim_{p\to
p_0,\, p\in\C_I} (p-p_0)^{-1}(G_\Lambda(p)-G_{\Lambda}(p_0))
$$
exist. Thus the functions $F_\Lambda$ and $G_\Lambda$ are
holomorphic on $\Omega\cap\C_I$ and so they admit a Cauchy formula on
the plane $\C_I$, computed e.g. on a circle $\gamma$, contained
in $\C_I$, whose interior contains $p_0$ and is contained in
$\Omega$.  Note that $(\Lambda f)_I(x+Iy)=\Lambda f_I(x+Iy)$. Moreover, if $p_0$ is real we can pick any complex plane
$\C_I$.  For any increment $h$ in $\C_I$ we compute
\[
\begin{split}
&h^{-1}((\Lambda f)_I(p_0+h)- (\Lambda f)_I(p_0))- \partial_s (\Lambda f)_I(p_0))
\\
&=\frac{1}{2\pi} \int_\gamma \left[h^{-1}
\left(\frac{1}{p-(p_0+h)}- \frac{1}{p-p_0}\right)-
\frac{1}{(p-p_0)^2}\right] dp_I (\Lambda f)_I(p)),\\
\end{split}
\]
where $dp_I=(dx+Idy)/I$.
Then we observe that $(\Lambda f)_I(p)$ is continuous on $\gamma$
which is  compact, so $|(\Lambda f)_I(p)|\leq C_\Lambda$ for all
$p\in\gamma$. The family of maps $f(p): \mathcal{X}^*\to \H$ is pointwise
bounded at each $\Lambda$, thus $\sup_{p\in\gamma}\|f_i (p)\|\leq C$
by the uniform boundedness theorem, see \cite{acs4}. Thus we have
$$
\left|\Lambda(h^{-1}(f_I (p_0+h)-f_I (p_0))- \partial_s (\Lambda f)_I(p_0)\right|
$$
$$
\leq\frac{C}{2\pi} \|\Lambda\| \int_\gamma \left|
\left(\frac{1}{p-(p_0+h)}- \frac{1}{p-p_0}\right)-
\frac{1}{(p-p_0)^2}\right| dp_I,
$$
so $h^{-1}(f_I(p_0+h)-f_I(p_0))$ is uniformly Cauchy for
$\|\Lambda\|\leq 1$ and by Lemma \ref{techlemma} it converges in $\mathcal{X}$.
Thus $f$ admits slice derivative at every $p_0\in \Omega$ and so it is
strongly slice hyperholomorphic in $\Omega$.
\end{proof}
\begin{Dn}
Let $\mathcal X$ be a two-sided Banach space over
$\mathbb H$.
We say that a function $f:\, \Omega\to \mathcal X$ is (weakly) slice
hypermeromorphic if for any $\Lambda\in \mathcal X^*$ the function $\Lambda f: \, \Omega\to
\mathbb H$ is slice hypermeromorphic in $\Omega$.
\end{Dn}
\begin{Rk}{\rm
The previous definition means, in particular, that $f:\, \Omega' \to \mathcal X$ is slice hyperholomorphic,
where the points in $\Omega\setminus \Omega'$ are the poles of $f$ and $\Omega\setminus \Omega'$ has no point limit in $\Omega$. }
\end{Rk}

Our next task is to prove that weakly slice hyperholomorphic
functions are those functions whose restrictions to any complex
plane $\C_I$ are in the kernel of the Cauchy-Riemann operator
$\partial_x +I\partial_y$.
\begin{Pn}\label{CRoperator}
Let $\mathcal X$ be a two sided quaternionic Banach space.
A real differentiable function $f: \Omega\subseteq\H\to \mathcal{X}$ is weakly
slice hyperholomorphic if and only if $(\partial_x
+I\partial_y)f_I(x+Iy)=0$ for all $I\in\mathbb{S}$.
\end{Pn}
\begin{proof}
If $f$ is weakly slice hyperholomorphic, then, as it happens in
the classical complex case, for every nonreal $p_0\in \Omega$,
$p_0\in\C_I$, we can compute the limit (\ref{deriv}) for the
function $\Lambda f_I$ choosing $p=p_0+h$ with $h\in\R$ and for
$p=p_0+Ih$  with $h\in\R$. We obtain, respectively, $\partial_x
f_I \Lambda (p_0)$ and $-I\partial_y \Lambda f_I(p_0)$ which
coincide. Thus we get $(\partial_x+I\partial_y)\Lambda
f_I(p_0)=\Lambda (\partial_x+I\partial_y) f_I(p_0)=0$ for any $\Lambda\in
\mathcal{X}^*$ and the statement follows by the Hahn-Banach theorem. If
$p_0$ is real, then the statement follows by an analogous
argument since the limit (\ref{deriv1}) exists for all
$I\in\mathbb{S}$.  Conversely, if $f_I$ satisfies the
Cauchy-Riemann on $\Omega\cap\C_I$ then $\Lambda((\partial_x
+I\partial_y)f_I(x+Iy))=0$ for all $\Lambda\in \mathcal{X}^*$ and all
$I\in\mathbb{S}$. Since $\Lambda$ is linear and continuous  we
can write $(\partial_x +I\partial_y)\Lambda f_I(x+Iy)=0$ and thus
the function $\Lambda f_I(x+Iy)$ is in the kernel of
$\partial_x+I\partial_y$ for all $\Lambda\in \mathcal{X}^*$ or,
equivalently by Proposition \ref{equivalence}, it admits slice
derivative. Thus at every $p_0\in \Omega\cap C_I$  we have
$$
\lim_{p \to p_0, p\in\C_I} (p-p_0)^{-1}( \Lambda f_I(p)- \Lambda
f_I(p_0))= \lim_{p \to p_0, p\in\C_I}  \Lambda( (p-p_0)^{-1}(
f_I(p)-f_I(p_0))),
$$
for all $\Lambda\in \mathcal{X}^*$. It follows that $f$ is  weakly slice
hyperholomorphic.
\end{proof}
Since the class of weakly and strongly slice hyperholomorphic
 functions coincide and in view of Proposition \ref{CRoperator}, from now on we will refer to them simply as slice
  hyperholomorphic functions and we denote the set of $\mathcal X$-valued slice hyperholomorphic functions on $\Omega$ by $\mathscr{S} (\Omega, \mathcal X)$.\\
 The following result is immediate:
  \begin{Pn}
Let $\mathcal{X}$ be a two sided quaternionic Banach space. Then the  set of slice hyperholomorphic functions defined on $\Omega\subseteq \mathbb{H}$ with values
in $\mathcal{X}$ is a right quaternionic linear space.
\end{Pn}
\begin{Pn}[Identity Principle] Let $\mathcal{X}$ be a two sided quaternionic Banach space,  $\Omega$ be an s-domain and let
$f,g: \Omega\subseteq\H \to \mathcal{X}$ be two slice hyperholomorphic functions.
 If $f=g$ on a set $Z\subseteq \Omega\cap \C_I$ having an accumulation point, for some $I\in\mathbb S$, then $f=g$ on $\Omega$.
\end{Pn}
\begin{proof}
The hypothesis implies $\Lambda f=\Lambda g$ on $Z$ for every
$\Lambda\in \mathcal{X}^*$ thus the slice hyperholomorphic function
$\Lambda (f-g)$ is identically zero not only on $Z$ but also on
$\Omega$ by the Identity Principle for quaternionic valued slice hyperholomorphic functions.
By the Hahn-Banach theorem $f-g=0$ on $\Omega$.
\end{proof}
\begin{Rk}{\rm
The Identity Principle implies that two slice hyperholomorphic functions defined on an s-domain and with values in a two sided quaternionic Banach space $\mathcal X$ coincide if their restrictions to the real axis coincide. More in general, any real analytic function $f: [a,b]\subseteq\mathbb R \to \mathcal X$ can be extended to a function ${\rm ext}(f)$ slice hyperholomorphic on an axially symmetric s-domain $\Omega$ containing $[a,b]$. The existence of the extension is assured by the fact that for any $x_0\in [a,b]$ the function $f$ can be written as $f(x)=\sum_{n\geq 0} x^n A_n$, $A_n\in\mathcal X$,  and $x$ such that $|x-x_0|<\varepsilon$ and thus $({\rm ext}f)(p)= \sum_{n\geq 0} p^n A_n$ for $|p-x_0|< \varepsilon_{x_0}$. Thus the claim holds setting $B(x_0,\varepsilon_{x_0})=\{p\in\mathbb H\ :\ |p-x_0|< \varepsilon_{x_0}\}$ and  $\Omega=\cup_{x_0\in I} B(x_0,\varepsilon_{x_0})$.
}
\end{Rk}

Let us recall, see \cite{MR2752913}, that the Cauchy kernel to be used in the Cauchy formula for slice hyperholomorphic functions is
$$
 S_L^{-1}(s,p)=-(p^2-2p{\rm Re}(s)+|s|^2)^{-1}(p-\overline{s}).
$$
It is a function slice hyperholomorphic on the left in the variable $p$ and on the right in $s$.
In the case of right regular functions the kernel  is
$$
S_R^{-1}(s,q):=-(q-\bar s)(q^2-2{\rm Re}(s)q+|s|^2)^{-1},
$$
which is slice hyperholomorphic on the right in the variable $q$ and on the left in $s$.
The Cauchy formula holds  for slice hyperholomorphic functions with values in a quaternionic Banach space:
\begin{Tm}[Cauchy formulas] \label{CauchygeneraleVV}
Let $\mathcal X$ be a  two sided quaternionic Banach space and let $W$  be an open set in $\hh$.
 Let $\overline{\Omega}\subset W$ be an axially symmetric s-domain,
 and let $\partial (\Omega\cap \mathbb{C}_I)$ be the union
of a finite number of rectifiable Jordan curves  for every
$I\in\mathbb{S}$. Set  $ds_I=ds/ I$.
If  $f:W \to \mathcal{X}$ is a left slice hyperholomorphic, then, for $q\in \Omega$, we have
\begin{equation}\label{cauchynuovo}
 f(p)=\frac{1}{2 \pi}\int_{\partial (\Omega\cap \mathbb{C}_I)} S_L^{-1}(s,p)ds_I f(s),
\end{equation}
if  $f:W \to \mathcal{X}$ is a right slice hyperholomorphic, then, for $q\in \Omega$, we have
\begin{equation}\label{Cauchyright}
 f(q)=\frac{1}{2 \pi}\int_{\partial (\Omega\cap \mathcal{C}_I)}  f(s)ds_I {S}_R^{-1}(s,q),
\end{equation}
and the integrals (\ref{cauchynuovo}), (\ref{Cauchyright})  do not depend on the choice of the imaginary unit $I\in\mathbb{S}$ and on  $\Omega\subset W$.
\end{Tm}
\begin{proof}
We have proved that weakly slice hyperholomorphic functions are strongly slice hyperholomorphic functions, so in particular they are continuous functions, so the validity of the formulas (\ref{cauchynuovo}), (\ref{Cauchyright}) follows as in point (b) p. 80 \cite{rudin}.
\end{proof}

We now show another description  of the class $\mathscr{S}(\Omega,\mathcal{X})$ of
slice hyperholomorphic functions on $\Omega$  with values in $\mathcal{X}$.
\begin{Dn}
Consider the set of functions of the form
$f(p)=f(x+Iy)=\alpha (x,y) +I\beta (x,y)$ where $\alpha, \beta:
\Omega\to \mathcal{X}$ depend only on $x,y$, are real differentiable, satisfy the Cauchy-Riemann
equations $\partial_x \alpha -\partial_y\beta=0$, $\partial_y \alpha
+\partial_x\beta=0$ and assume that $\alpha(x,-y)=\alpha(x,y)$,
$\beta(x,-y)=-\beta(x,y)$.
We will denote the class of function of this form by
$\mathscr{H}(\Omega,\mathcal{X})$.
\end{Dn}
Observe that the conditions on $\alpha$ and $\beta$ are required in order to have that the function $f$ is well posed. Note also that if $p=x$ is a real
quaternion, then $I$ is not uniquely defined but the hypothesis
that $\beta$ is odd in the variable $y$ implies $\beta(x,0)=0$.
\begin{Tm} Let $\Omega$ be an axially symmetric s-domain and let $\mathcal X$ be a two sided quaternionic Banach space. Then $\mathscr{S}(\Omega,\mathcal{X})=\mathscr{H}(\Omega,\mathcal{X})$.
\end{Tm}
\begin{proof} The inclusion  $\mathscr{H}(\Omega,\mathcal{X})\subseteq \mathscr{S}(\Omega,\mathcal{X})$ is clear: any function $f\in \mathscr{H}(\Omega,\mathcal{X})$
is real differentiable and such that $f_I$ satisfies $(\partial_x+I\partial_y)f_I=0$ (note that this implication
does not need any hypothesis on the open set $\Omega$).
Conversely, assume that $f\in\mathscr{S}(\Omega,\mathcal{X})$. Let us show that
$$
f(x+Iy)=\frac 12 (1-IJ)f(x+Jy)+ \frac 12 (1+IJ)f(x-Jy).
$$
If we consider real quaternions, i.e. $y=0$, then the formula holds trivially.
For nonreal quaternions, set
$$
\phi(x+Iy)=\frac 12 (1-IJ)f(x+Jy)+ \frac 12 (1+IJ)f(x-Jy).
$$
Then, using the fact that $f$ is slice hyperholomorphic, it is
immediate that $(\partial_x+I\partial y)\phi(x+Iy)=0$ and so $\phi$ is slice hyperholomorphic. Since
$\phi=f$ on $\Omega\cap\C_I$ then it coincides with $f$ on $\Omega$ by the
Identity Principle. By writing
$$
f(x+Iy)=\frac 12 \left[(f(x+Jy)+f(x-Jy)+  IJ
(f(x-Jy)-f(x+Jy))\right]
$$
 and setting $\alpha(x,y)=\frac 12 (f(x+Jy)+f(x-Jy))$,
 $\beta(x,y)=\frac 12 J (f(x-Jy)-f(x+Jy))$ we have that $f(x+Iy)=\alpha(x,y)+I\beta(x,y)$.
 Reasoning as in \cite[Theorem 2.2.18]{MR2752913}
 we see that $\alpha$, $\beta$ do not depend on $I$.
 It is then an easy computation
 to verify that $\alpha$, $\beta$ satisfy the above assumptions.
\end{proof}
Using this alternative description of slice hyperholomorphic
functions with values in $\mathcal{X}$, we can now define a notion
of product which is based on a suitable pointwise multiplication. To this end we need
an additional structure on the two sided quaternionic Banach space $\mathcal{X}$. Assume that in $\mathcal X$ is defined
a multiplication which is associative, distributive with respect to the sum in $\mathcal X$. Assume also that
$q(x_1x_2)=(qx_1)x_2$ and $(x_1x_2)q=x_1(x_2q)$ for all $q\in\mathbb H$ and for all $x_1,x_2\in\mathcal X$. Then we will say that $\mathcal X$
is a two sided quaternionic Banach algebra. As customary we will say that the algebra $\mathcal X$ is with unity is $\mathcal X$ possesses a unity with respect to the product.
\begin{Dn}
Let $\Omega\subseteq\H$ be an axially symmetric s-domain and let
$f,g:  \Omega\to \mathcal{X}$ be slice hyperholomorphic functions with values in a two sided quaternionic Banach algebra $\mathcal X$. Let
$f(x+Iy)=\alpha(x,y)+I\beta(x,y)$,
$g(x+Iy)=\gamma(x,y)+I\delta(x,y)$. Then we define
\begin{equation}\label{starleft}
(f\star g)(x+Iy):= (\alpha\gamma -\beta \delta)(x,y)+
I(\alpha\delta +\beta \gamma )(x,y).
\end{equation}
\end{Dn}
By construction, the function $f\star g$ is slice
hyperholomorphic, as it can be easily verified.
\begin{Rk}{\rm
If $\Omega$ is a ball with center at a real point (let us assume at
the origin for simplicity) then it  is immediate that $f$, $g$ admit power series
expansion  and thus if $f(p)=\sum_{n= 0}^\infty
p^n a_n$, $g(p)=\sum_{n= 0}^\infty p^n b_n$, $a_n,b_n\in \mathcal{X}$ for all
$n$. Then $f\star g(p):=\sum_{n= 0}^\infty p^n (\sum_{r=0}^n
a_rb_{n-r})$ where the series converges. }
\end{Rk}
\begin{Rk}{\rm In case we consider right slice hyperholomorphic functions, the class $\mathscr{H}(\Omega,\mathcal{X})$ consists of functions of the form $f(x+Iy)=\alpha(x,y)+\beta (x,y)I$ where $\alpha, \beta$ satisfy the assumptions discussed above. We now give the notion of right slice product, denoted by $\star_r$. Given two right slice hyperholomorphic functions
$f,g:  \Omega\to \mathcal{X}$ with values in a two sided quaternionic Banach algebra $\mathcal X$ where
$f(x+Iy)=\alpha(x,y)+\beta(x,y)I$,
$g(x+Iy)=\gamma(x,y)+\delta(x,y)I$, we define
\begin{equation}\label{starright}
(f\star_r g)(x+Iy):= (\alpha\gamma -\beta \delta)(x,y)+
(\alpha\delta +\beta \gamma )(x,y) I.
\end{equation}
 }
\end{Rk}
\begin{Rk}
{\rm It is important to point out that if one is in need of considering slice hyperholomorphic functions on axially symmetric open sets $U$ which are not necessarily s-domains, then it is more convenient to use the class $\mathcal{H}(\Omega,\mathcal X)$ because they allow to have a notion of multiplication.}
\end{Rk}

\begin{Rk}{\rm
Consider the following case:
let $\Omega$ be an axially symmetric s-domain in $\mathbb H$ and
let $\mathcal H_i$; $i=1,2,3$ be two sided quaternionic Hilbert
spaces. Let $f:\, \Omega \to \mathbf L(\mathcal H_1,\mathcal
H_2)$, $g:\, \Omega \to \mathbf L(\mathcal H_2,\mathcal H_3)$ be
slice hyperholomorphic and let
$$
f(p)=f(x+Iy)=\alpha(x,y)+I\beta(x,y),\ \ \
g(p)=g(x+Iy)=\gamma(x,y)+I\delta(x,y).
$$
 We define the $\star$-product as in \eqref{starleft}
If $f,g$ are right slice hyperholomorphic, then
 we define the $\star_r$-product as in \eqref{starright}.
The product $\alpha(x,y)\gamma(x,y)$  (and the other three
products appearing in $f\star g$) is an operator belonging to
$\mathbf L(\mathcal H_1,\mathcal H_3)$, thus
 $f\star g: \, \Omega \to \mathbf L(\mathcal H_1,\mathcal H_3)$.
In the special case in which
$$
f(p)=\sum_{n=0}^\infty p^n A_n, \ \ \ \ A_n\in \mathbf L(\mathcal H_1,\mathcal H_2),
$$
$$g(p)=\sum_{n=0}^\infty p^n B_n, \ \ \ \ B_n\in \mathbf L(\mathcal H_2,\mathcal H_3),
$$
 then we have
 $$f\star g(p)=\sum_{n=0}^\infty p^n(\sum_{r=0}^n A_rB_{n-r}),$$ as expected.}
\end{Rk}

\subsection{Realizations}
The following notions of S-spectrum and of  S-resolvent operator  will be used in the sequel.
\begin{Dn}
Let $\mathcal X$ be a two sided quaternionic Banach space and let
$A$ be a bounded operator on $\mathcal X$ into itself.
We define the $S$-spectrum $\sigma_S(A)$ of $A$  as:
$$
\sigma_S(A)=\{ p\in \mathbb{H}\ \ :\ \ A^2-2 \ {\rm
Re}\,(p)A+|p|^2 I\ \ \ {\it is\ not\  invertible}\}.
$$
The $S$-resolvent set $\rho_S(A)$ is defined by $\rho_S(A)=\mathbb{H}\setminus
\sigma_S(A).$
\\

For $p\in \rho_S(A)$  the right $S$-resolvent operator is defined
as
\begin{equation}\label{quatSresorig}
S_R^{-1}(p,A):=-(A-\overline{p}I)(A^2-2{\rm Re}\,(p) A+|p|^2
I)^{-1}.
\end{equation}
\end{Dn}

\begin{Rk}\label{spettro1}{\rm It is useful to recall that when $A$ is a matrix its (point) $S$-spectrum coincides with its right spectrum, see e.g.
\cite{acs3}. When $p\in\mathbb R$ or, more in general, when $p$ commute with an operator $A$, then $S^{-1}_R(p,A)=(pI-A)^{-1}$, see Proposition 3.1.6 in \cite{MR2752913}. }
\end{Rk}

\begin{Pn}\label{extA}
 Let $\mathcal X$ be a two sided quaternionic Banach space and let $f:\,\rho_S(A)\cap \mathbb R\setminus\{0\}\to \mathcal X$ be the function $f(x)=(I-xA)^{-1}$.
Then
\[
p^{-1}S_R^{-1}(p^{-1},A)=(I-\bar p A)(I-2 {\rm Re}(p) A
+|p|^2 A^2)^{-1}
\]
is the unique slice hyperholomorphic extension to $\rho_S(A)$.
\end{Pn}
\begin{proof}
The fact that $p^{-1}S_L^{-1}(p^{-1},A)$ is slice hyperholomorphic
in $p$ outside the S-spectrum is trivial since it is the
S-resolvent and it coincides with the function $f$ on the real
axis. The uniqueness follows from the identity principle.
\end{proof}

The notation $S^{-1}_R(p^{-1},A)$ comes
from \cite{MR2752913} but we will also write
\[
p^{-1}S^{-1}_R(p^{-1},A)=(I-pA)^{-\star}.
\]

This last expression makes sense when $A$ acts on a two-sided
quaternionic vector space. In a more general setting, we have the
following result:
\begin{Pn}
Let $A$ be a bounded linear operator from a right-sided
quaternionic Banach $\mathcal P$ space into itself, and let $G$
be a bounded linear operator from $\mathcal P$ into $\mathcal Q$,
where $\mathcal Q$ is a  two sided quaternionic Banach space. The
slice hyperholomorphic extension of $G(I-xA)^{-1}$,
$1/x\in\sigma_S(A)\cap\mathbb{R}$, is
\[
(G-\overline{p}GA)(I-2{\rm Re}(p)\, A +|p|^2A^2)^{-1}.
\]
\label{formula060813}
\end{Pn}
\begin{proof}
First we observe that $G(I-xA)^{-1}=\sum_{n= 0}^\infty x^nGA^n$
for $|x|\|A\|<1$. It is immediate that, for $|p|\|A\|<1$, the slice
hyperholomorphic extension of the series $\sum_{n= 0}^\infty
x^nGA^n$ is $\sum_{n= 0}^\infty p^nGA^n$ (as it is a converging power series
with coefficients on the right). To show that
\[
\sum_{n= 0}^\infty p^nGA^n=(G-\overline{p}GA)(I-2{\rm Re}(p)\, A
+|p|^2A^2)^{-1}
\]
we prove instead the equality
\[
(\sum_{n= 0}^\infty p^nGA^n)(I-2{\rm Re}(p)\, A+|p|^2A^2)=(G-\overline{p}GA).
\]
The left hand side gives
\[
\begin{split}
&\sum_{n= 0}^\infty p^nGA^n -2 \sum_{n= 0}^\infty {\rm Re}(p) p^nGA^{n+1}+ \sum_{n= 0}^\infty |p|^2p^nGA^{n+2}
\\
&= G + (p -2{\rm Re}(p)) GA +(p^2-2p{\rm Re}(p) +|p|^2)\sum_{n=0
}^\infty p^n GA^{n+2}\\
&= G -\bar p GA
\end{split}
\]
where we have used the identity $p^2-2p{\rm Re}(p) +|p|^2=0$. This completes the proof.
\end{proof}
\begin{Rk}{\rm
In analogy with the matrix case we will write, with an abuse of
notation in this case,  $G\star (I-pA)^{-\star}$ instead of the
expression $(G-\overline{p}GA)(I-2{\rm Re}(p)\, A
+|p|^2A^2)^{-1}$. } \label{milano2013}
\end{Rk}

\begin{Pn}
With the notation in Remark \ref{milano2013} it holds that
\begin{equation}
\label{formreal1}
D+pC\star(I-pA)^{-1}
B=D^{-1}-pD^{-1}C\star(I-p(A-BD^{-1}C))^{-\star}BD^{-1},
\end{equation}
and
\begin{equation}
\label{formreal2}
\begin{split}
(D_1+pC_1\star(I-pA_1)^{-\star}B_1)\star(D_2+pC_2\star(I-pA_2)^{-\star}B_2)&=\\
\hspace{-5cm}=D_1D_2+p\begin{pmatrix}C_1&D_1C_2\end{pmatrix}\star
\left(I-p\begin{pmatrix}A_1&B_1C_2\\0&A_2\end{pmatrix}\right)^{-\star}
\begin{pmatrix}B_1D_2\\ B_2\end{pmatrix}.
\end{split}
\end{equation}
\end{Pn}
\begin{proof}
When $p$ is real, the $\star$-product is replaced by the operator
product (or matrix product in the finite dimensional case) and
formulas \eqref{formreal1} and \eqref{formreal2} are then well
known, see e.g. \cite{bgr1}. Slice-hyperholomorphic extensions
lead then to the required result.
\end{proof}

\section{The Hardy space of the half-space $\mathbb H_+$}
\setcounter{equation}{0}

Let $\Pi_+$ be the right open half-plane  of
complex numbers $z$ such that ${\rm Re}(z)>0$.
The Hardy space  $\mathbf{H}_2(\Pi_+)$ consists of functions $f$ holomorphic in $\Pi_+$ such that
\begin{equation}\label{Hardy1}
\sup_{x>0}\int_{-\infty}^\infty |f(x+iy)|^2 dy <\infty.
\end{equation}
We recall that a function $f\in \mathbf{H}_2(\Pi_+)$ has nontangential limit $f(iy)$ for almost all $iy$ on the imaginary axis and
$f(iy)\in L_2(\mathbb R)$, see \cite[Theorem 3.1]{garnett}, moreover
\begin{equation}\label{Hardy2}
\sup_{x>0}\int_{-\infty}^\infty |f(x+iy)|^2 dy = \int_{-\infty}^\infty |f(iy)|^2 dy .
\end{equation}
Let us consider the kernel
\[
k_{\Pi_+}(z,w)=\frac{1}{2\pi}\frac{1}{z+\bar w},
\]
which is positive definite on $\Pi_+$. Then, the associated
reproducing kernel Hilbert space
   is the Hardy space  $\mathbf{H}_2(\Pi_+)$ endowed with
the scalar product
\[
\langle f,g\rangle_{\mathbf{H}_2(\Pi_+)}=\int_{-\infty}^{+\infty}
\overline{g(iy)}  f(iy)  dy,
\]
where $f,g\in \mathbf{H}_2(\Pi_+)$,
and the norm in  $\mathbf{H}_2(\Pi_+)$ is given by
$$
\|f\|_{\mathbf{H}_2(\Pi_+)}=\left( \int_{-\infty}^{+\infty} |
f(iy) |^2 dy  \right)^{\frac 12} .
$$
The kernel $k_{\Pi_+}(z,w)$ is reproducing in the sense that for every
$f\in \mathbf{H}_2(\Pi_+)$
\[
f(w)=\langle f(z),k_{\Pi_+}(z,w)\rangle_{\mathbf{H}_2(\Pi_+)}=\int_{-\infty}^\infty k_{\Pi_+}(w,iy) f(iy) dy,
\]
Let us now consider the half-space $\mathbb H_+$ of the quaternions
$q$ such that ${\rm Re}(q)>0$ and set $\Pi_{+,I}=\mathbb H_+\cap
\mathbb{C}_I$. We will denote by $f_I$ the restriction of a function $f$ defined on $\mathbb H_+$ to $\Pi_{+,I}$. We define
\[
\mathbf{H}_2(\Pi_{+,I})=\{f\ {\rm slice\ hyperholomorphic\ in\ \mathbb H_+}\ : \
\int_{-\infty}^{+\infty}   |f_I(Iy)|^2 dy <\infty\},
\]
where $f(Iy)$ denotes the nontangential value of $f$ at $Iy$. Note that these value exist almost everywhere, in fact any  $f\in \mathbf{H}_2(\Pi_{+,I})$ when restricted to a complex plane $\mathbb C_I$ can be written as $f_I(x+Iy)=F(x+Iy)+G(x+Iy)J$ where $J$ is any element in $\mathbb S$ orthogonal to $I$, and $F,G$ are $\mathbb C_I$-valued holomorphic functions.
Since the nontangential values of $F$ and $G$ exist almost everywhere at $Iy$, also the nontangential value of $f$ exists at $Iy$ a. e. on $\Pi_{+,I}$ and $f_I(Iy)=F(Iy)+G(Iy)J$ a.e.
\begin{Rk}{\rm In alternative, we could have defined  $\mathbf{H}_2(\Pi_{+,I})$ as the set of slice hyperholomorphic functions $f$ such that
$\sup_{x>0}\int_{-\infty}^{+\infty}   |f_I(x+Iy)|^2 dy <\infty$. However note that  $f_I(x+Iy)=F(x+Iy)+G(x+Iy)J$, see the above discussion, and so $|f_I(x+Iy)|^2=|F(x+Iy)|^2+|G(x+Iy)|^2$. Thus, using \eqref{Hardy2}, we have
\begin{equation}\label{HardyH}
\begin{split}
\sup_{x>0}\int_{-\infty}^{+\infty}   |f_I(x+Iy)|^2 dy&=\sup_{x>0}\int_{-\infty}^{+\infty}   |F(x+Iy)|^2 dy+\sup_{x>0}\int_{-\infty}^{+\infty}   |G(x+Iy)|^2 dy\\
&=\int_{-\infty}^{+\infty}   |F(Iy)|^2 dy+\int_{-\infty}^{+\infty}   |G(Iy)|^2 dy\\
&=\int_{-\infty}^{+\infty}   |f_I(Iy)|^2 dy.
\end{split}
\end{equation}
}
\end{Rk}
In $\mathbf{H}_2(\Pi_{+,I})$ we define the scalar product
\[
\langle
f,g\rangle_{\mathbf{H}_2(\Pi_{+,I})}=\int_{-\infty}^{+\infty}
\overline{g_I(Iy)} f_I(Iy) dy ,
\]
where $f_I(Iy)$, $g_I(Iy)$ denote the nontangential values of $f,g$
at $Iy$ on $\Pi_{+,I}$. This scalar product gives the norm
\[
\|f\|_{\mathbf{H}_2(\Pi_{+,I})}=\left(\int_{-\infty}^{+\infty}  |
f_I(Iy)|^2 dy\right)^{\frac 12},
\]
(which is finite by our assumptions).
\begin{Pn} Let $f$ be slice hyperholomorphic in $\mathbb H_+$ and assume that $f\in\mathbf{H}_2(\Pi_{+,I})$ for some $I\in\mathbb S$. Then
 for all
$J\in\mathbb S$ the following inequalities hold
$$
\frac{1}{2}\|f\|_{\mathbf{H}_2(\Pi_{+,I})}\leq  \|f\|_{\mathbf{H}_2(\Pi_{+,J})} \leq 2\|f\|_{\mathbf{H}_2(\Pi_{+,I})}.
$$
\end{Pn}
\begin{proof}
Formula (\ref{repr}) implies the inequality
$$
|f(x+Jy)|\leq |f(x+Iy)|+|f(x-Iy)|,
$$
and also
\begin{equation}\label{reprsquared}
|f(x+Jy)|^2\leq 2(|f(x+Iy)|^2+|f(x-Iy)|^2).
\end{equation}
Using \eqref{HardyH}, \eqref{repr} and \eqref{reprsquared} we deduce
\[
\begin{split}
\|f
\|_{\mathbf{H}_2(\Pi_{+,J})}^2=\int_{-\infty}^{+\infty}  |
f_J(Jy)|^2 dy&=\sup_{x>0} \int_{-\infty}^{+\infty}  |
f_J(x+Jy)|^2 dy\\
&\leq \sup_{x>0} \int_{-\infty}^{+\infty}
2(|f_I(x+Iy)|^2+f_I(x-Iy)|^2) dy\\
&=4\int_{-\infty}^{+\infty}
|f_I(Iy)|^2 dy
\end{split}
\]
and so
 $\|f\|^2_{\mathbf{H}_2(\Pi_{+,J})}\leq 4
\|f\|^2_{\mathbf{H}_2(\Pi_{+,I})}$. By changing $J$ with $I$ we
obtain the reverse inequality
and the statement follows.
\end{proof}
An immediate consequence of this result is:
\begin{Cy} A function $f\in \mathbf{H}_2(\Pi_{+,I})$ for some $I\in\mathbb S$ if and only if $f\in \mathbf{H}_2(\Pi_{+,J})$ for all $J\in\mathbb S$.
\end{Cy}
We now introduce the Hardy space of the half space $\mathbb H_+$:
\begin{Dn}
We define $\mathbf{H}_2(\mathbb H_+)$ as the space of slice
hyperholomorphic functions on $\mathbb H_+$ such that
\begin{equation}\label{xge0integ}
\sup_{I\in\mathbb S}\int_{-\infty}^{+\infty} | f(Iy) |^2 dy  <\infty .
\end{equation}
\end{Dn}
We have:
\begin{Pn}
The function
\begin{equation}\label{kernel}
k(p,q)=(\bar p +\bar q)(|p|^2 +2{\rm Re}(p) \bar q
+\bar q^2)^{-1}
\end{equation}
is slice hyperholomorphic in $p$ and $\bar q$ on the left and on
the right, respectively in its domain of definition, i.e. for $p\not\in[\bar q]$.
The restriction of $\frac{1}{2\pi}k(p,q)$ to $\mathbb{C}_I\times
\mathbb{C}_I$ coincides with $k_{\Pi_+}(z,w)$. Moreover we have the
equality:
\begin{equation}\label{kernel1}
k(p,q)=(|q|^2 +2{\rm Re}(q) p + p^2)^{-1}( p + q).
\end{equation}
\end{Pn}
\begin{proof}
Some computations allow to obtain $k(p,q)$ as the left slice hyperholomorphic
extension in $z$ of $k_q(z)=k(z,q)$, by taking $z$ on the same complex plane as $q$.
The function we obtain turns out to be also right slice hyperholomorphic in $\bar q$.
The second equality follows by taking the right slice hyperholomorphic extension in
$\bar q$ and observing that it is left slice hyperholomorphic in $p$.
\end{proof}
\begin{Pn}
The kernel $\frac{1}{2\pi}k(p,q)$ is reproducing, i.e. for any
$f\in\mathbf{H}_2(\mathbb H_+)$
\[
f(p)=\int_{-\infty}^\infty \frac{1}{2\pi}k(p,Iy) f(Iy) dy.
\]
\end{Pn}
\begin{proof}
Let $q=u+ I_q v$ and let $p=u+Iv$ be the point on the sphere
determined by $q$ and belonging to the plane $\mathbb C_I$. Then
we have
\[
f(p)=\int_{-\infty}^\infty \frac{1}{2\pi}k(p,Iy) f(Iy) dy, \qquad f(\bar
p)=\int_{-\infty}^\infty \frac{1}{2\pi}k(\bar p,Iy) f(Iy) dy.
\]
The extension formula \eqref{ext} applied to $k_{Iy}(p)=k(p,Iy)$ shows the
statement.
\end{proof}
The following property will be useful in the sequel:
\begin{Pn}
\label{identity}
The kernel $k(p,q)$ satisfies
$$
pk(p,q)+k(p,q)\overline{q}=1.
$$
\end{Pn}
\begin{proof}
From the expression (\ref{kernel}), and since $q$ commutes with
$(|p|^2 +2{\rm Re}(p) \bar q +\bar q^2)^{-1}$, we have
\[
\begin{split}
&p(\bar p +\bar q)(|p|^2 +2{\rm Re}(p) \bar q +\bar q^2)^{-1}+(\bar p +\bar q)(|p|^2 +2{\rm Re}(p) \bar q +\bar q^2)^{-1}\overline{q}\\
&=(|p|^2+p\bar q +\bar p\bar q +\bar q^2)(|p|^2 +2{\rm Re}(p)
\bar q +\bar q^2)^{-1}=1.
\end{split}
\]
\end{proof}
We know that if $\{\phi_n(z)\}$ is an orthonormal basis for
$\mathbf{H}_2(\Pi_{+,I})$, for some $I\in\mathbb S$, then
\begin{equation}\label{kappa}
k(z,w)=\sum_{n=1}^\infty \phi_n(z) \overline{\phi_n(w)},
\end{equation}
and so the kernel $k(z,w)$ is positive definite. We now prove the following:
\begin{Pn}
Let $\{\phi_n(z)\}$ be an orthonormal basis for
$\mathbf{H}_2(\Pi_{+,I})$, for some $I\in\mathbb S$,  and let
$\{\Phi_n(q)\}=\{{\rm ext}(\phi_n(z))\}$ be the sequence of  the
slice hyperholomorphic extensions of its elements. Then
$\{\Phi_n(q)\}$ is an orthonormal basis for
$\mathbf{H}_2(\mathbb H_+)$,  and
$$
k(p,q)=\sum_{n=1}^\infty \Phi_n(p)\overline{\Phi_n(q)}.
$$
\end{Pn}
\begin{proof}
Let $\{\phi_n(z)\}$ be an orthonormal basis for
$\mathbf{H}_2(\Pi_{+,I})$ and let $\{\Phi_n(q)\}=\{{\rm
ext}(\phi_n(z))\}$ be the sequence of  the slice hyperholomorphic
extensions of its elements. Then $\{\Phi_n(q)\}$ is a
generating set for $\mathbf{H}_2(\mathbb H_+)$. In fact take any
$f\in \mathbf{H}_2(\mathbb H_+)$ and consider its restriction to a
complex plane $\mathbb C_I$, for some $I\in\mathbb S$. Then, by
choosing $J\in\mathbb S$ such that $I,J$ are orthogonal, and taking $q=x+Iy$ we have
$f_I(x+Iy)= F(x+Iy)+G(x+Iy)J$
with $F,G$ holomorphic on $\Pi_{+,I}$ and
$$
\int_{-\infty}^{+\infty} | f(Iy) |^2 dy =\int_{-\infty}^{+\infty}
(| F(Iy) |^2 +| G(Iy) |^2) dy  <\infty
$$
and, as a consequence,
$$
\int_{-\infty}^{+\infty} | F(Iy) |^2 dy
\leq\int_{-\infty}^{+\infty} (| F(Iy) |^2 +| G(Iy) |^2) dy
<\infty.
$$
We deduce that both
$$
\int_{-\infty}^{+\infty} | F(Iy) |^2 dy \quad {\rm and}
\quad \int_{-\infty}^{+\infty} | G(Iy) |^2 dy
$$
are finite and so $F,G$ belong to $\mathbf{H}_2(\Pi_{+,I})$. We
can write $F(x+Iy)=\sum_{n=1}^\infty \phi_n(z) a_n$ and
$G(x+Iy)=\sum_{n=1}^\infty \phi_n(z) b_n$, thus
$f_I(x+Iy)= \sum_{n=1}^\infty \phi_n(z) (a_n+b_nJ)$. By taking
the extension with respect to $z$ we finally obtain $f(q)=
\sum_{n=1}^\infty \Phi_n(q) (a_n+b_n J)$. The fact that
$\{\Phi_n(p)\}$ is made by orthonormal elements (thus linearly
independent) in $\mathbf{H}_2(\mathbb H_+)$ follows from
\[
\begin{split}
\langle \Phi_n(p), \Phi_m(p)\rangle_{\mathbf{H}_2(\Pi_{+,I_p})}& =\int_{-\infty}^\infty \overline{\Phi_m(I_py)} \Phi_n(I_py) dy\\
&=\int_{-\infty}^\infty \overline{\phi_m(I_py)} \phi_n(I_py)
dy=\delta_{nm}.
\end{split}
\]
Then (\ref{kappa}) yields
$$
k(p,w)={\rm ext}_zk(z,w)=\sum_{n=1}^\infty {\rm
ext}_z(\phi_n)(z) \overline{\phi_n(w)}=\sum_{n=1}^\infty
\Phi_n(p) \overline{\phi_n(w)},
$$
where we have written ${\rm ext}_z$ to emphasize that we are
taking the extension in the variable $z$ (note that in this way
we have obtained the kernel written in the form (\ref{kernel})).
Now we observe that the function $\sum_{n=1}^\infty
\Phi_n(p) \overline{\Phi_n(q)}$ is slice hyperholomorphic on the left and on the right with respect to $p$ and $\bar q$, respectively, and coincides with $k(p,w)$ when restricted to the plane containing $w$. By the uniqueness of the extension we have $k(p,q)= \sum_{n=1}^\infty
\Phi_n(p) \overline{\Phi_n(q)}$,
  and the statement follows.
\end{proof}

\bigskip

\noindent
We now introduce the Blaschke factors in the half space $\mathbb H_+$.
\begin{Dn}
For $a\in\mathbb H_+$ set
\[
b_{a}(p)=(p+\bar a)^{-\star}\star (p-a).
\]
The function $b_a(p)$ is called Blaschke factor at $a$ in the half space $\mathbb H_+$.
\end{Dn}
\begin{Rk}{\rm
The function $b_a(p)$ is defined outside the sphere $[-a]$ as it can be easily seen by rewriting it as
$$
b_a(p)=(p^2+2{\rm Re}(a)p+|a|^2)^{-1}  (p+a)\star (p-a)= (p^2+2{\rm Re}(a)p+|a|^2)^{-1}  (p^2-a^2)
$$
and it has a zero for $p=a$. Note in fact that $p=-a$ is not a zero since it is a pole (of order 0).
When $a\in\mathbb R$ the function $b_a(p)=(p+ a)^{-1}(p-a)$ has a pole at $p=- a$ and a zero at $a$.
A Blaschke factor is slice hyperholomorphic where it is defined, by construction.}
\end{Rk}
We have the following result which characterizes the convergence of a Blaschke product. We denote by $\Pi^\star$  the $\star$-product:
\begin{Tm}\label{productblaschke}
Let $\{a_j\}\subset \mathbb{H}_+$, $j=1,2,\ldots$ be a sequence
of quaternions such that  $\sum_{j\geq 1} {\rm Re} (a_j) <
\infty$. Then the function
\begin{equation}\label{B_producthp}
B(p):=\Pi^\star_{j\geq 1}   (p+\bar a_j)^{-\star}\star (p-a_j),
\end{equation}
converges
uniformly on the compact subsets of $\mathbb{H}_+$.
\end{Tm}
\begin{proof}
We reason as in the proof of the corresponding result in the complex case (but see also the proof of Theorem 5.6 in \cite{acs2}).
We note that, see Remark 5.4 in \cite{acs2}, we can write
\begin{equation}\label{blaschkesimpler}
(p+\bar a_j)^{-\star}\star (p-a_j)=(\tilde p+\bar a_j)^{-1} (\tilde p-a_j)
\end{equation}
where $\tilde p= \lambda^c(p)^{-1} p \lambda^c(p)$ and $\lambda^c(p)=p+a_j$ (note that $\lambda^c(p)\not=0$ for $p\not\in [-a_j]$)  and so
\begin{equation}\label{formula3.8}
(p+\bar a_j)^{-\star}\star (p-a_j)=(\tilde p+\bar a_j)^{-1} (\tilde p-a_j)=1-2 {\rm Re} (a_j) (\tilde p+\bar a_j)^{-1} .
\end{equation}
 By taking the modulus of the right hand side of \eqref{B_producthp}, using \eqref{formula3.8}, and reasoning as in the complex case, we conclude that the Blaschke product converges if and only if
$\sum_{j=1}^\infty {\rm Re}(a_j) <\infty$.
\end{proof}
As in the unit disk case, we have two kinds of Blaschke factors. In fact, products of the form
\[
b_{a}(p)\star b_{\bar a}(p)=((p+\bar a)^{-\star}\star (p-a))\star ((p+ a)^{-\star}\star (p-\bar a))
\]
can be written as
\[
b_{a}(p)\star b_{\bar a}(p)=(p^2 +2 {\rm Re}(a)p +|a|^2)^{-1} (p^2 -2 {\rm Re}(a)p +|a|^2),
\]
and they admit the sphere $[a]$ as set of zeros. Note that slice regular functions which vanish at two different points belonging to the same sphere in reality  vanish on the whole sphere (see \cite[Corollary 4.3.7]{MR2752913}. Thus if we want to construct a Blaschke product vanishing at some prescribed points and spheres, it is convenient to introduce the following:
\begin{Dn}
For $a\in\mathbb H_+$ set
\[
b_{[a]}(p)=(p^2 +2 {\rm Re}(a)p +|a|^2)^{-1} (p^2 -2 {\rm Re}(a)p +|a|^2).
\]
The function $b_a(p)$ is called Blaschke factor at the sphere $[a]$ in the half space $\mathbb H_+$.
\end{Dn}
Note that the definition is well posed since it does not depend on the choice of the point $a$. As a consequence of Theorem \ref{productblaschke} we have:
\begin{Cy}\label{Cyproduct}
Let $\{c_j\}\subset \mathbb{H}_+$, $j=1,2,\ldots$ be a sequence
of quaternions such that  $\sum_{j\geq 1} {\rm Re} (c_j) <
\infty$. Then the function
\begin{equation}\label{B_product}
B(p):=\Pi_{j\geq 1}   (p^2 +2 {\rm Re}(c_j)p +|c_j|^2)^{-1} (p^2 -2 {\rm Re}(c_j)p +|c_j|^2),
\end{equation}
converges
uniformly on the compact subsets of $\mathbb{H}_+$.
\end{Cy}
\begin{proof}
It is sufficient to write  $B(p)=\prod_{j\geq 1}b_{[c_j]}(p)=\prod_{j\geq 1}b_{c_j}(p)\star b_{\bar c_j}(p)$ and to observe that $\sum_{j\geq 1} {\rm Re}(c_j) <\infty$ by hypothesis.
\end{proof}
To state the next result, we need to repeat the notion of multiplicity of a sphere of zeros and of a point which is an isolated zero.\\
We say that the {\em multiplicity of the spherical zero} $[c_j]$ of a function $Q(p)$ is $m_j$ if  $m_j$ is the maximum of the integers $m$ such that $(p^2+2{\rm Re}(c_j)p+|c_j|^2)^{m}$ divides ${Q}(p)$.\\
Let $\alpha_j\in\mathbb H\setminus \mathbb R$ and let
\begin{equation}\label{multz}
{Q}(p)=(p-\alpha_1)\star \ldots \star (p-\alpha_n)\star g(p)\quad \alpha_{j+1}\not=\bar{\alpha}_{j},\ \  j=1,\ldots, n-1, \ \ g(p)\not=0.
\end{equation}
We say that $\alpha_1$ is a zero of $Q$ of {\em multiplicity} $1$ if $\alpha_j\not\in [\alpha_1]$ for $j=2,\ldots, n$.\\
We say that $\alpha_1$ is a zero of $Q$ of {\em multiplicity} $n\geq 2$ if $\alpha_j\in [\alpha_1]$ for all $j=2,\ldots, n$.\\
If $\alpha_j\in\mathbb R$ we can repeat the notion of multiplicity of $\alpha_1$ where \eqref{multz} holds by removing the assumption $\alpha_{j+1}\not=\bar{\alpha}_{j}$. This definition coincides with the standard notion of multiplicity since, in this case, the $\star$-product reduces to the pointwise product. Note that if a function has a sphere of zeros at $[\alpha]$ with multiplicity $m$, at most one point on $[\alpha]$ can have higher multiplicity; in fact if there are two such points it means that the sphere $[\alpha]$ of zeros has higher multiplicity.

Thus we can prove the following:
\begin{Tm}
A Blaschke product having zeros at the set
 $$
 Z=\{(a_1,\mu_1), (a_2,\mu_2), \ldots, ([c_1],\nu_1), ([c_2],\nu_2), \ldots \}
 $$
 where $a_j\in \mathbb{H}_+$, $a_j$ have respective multiplicities $\mu_j\geq 1$, $[a_i]\not=[a_j]$ if $i\not=j$, $c_i\in \mathbb{H}_+$, the spheres $[c_j]$ have respective multiplicities $\nu_j\geq 1$,
 $j=1,2,\ldots$, $[c_i]\not=[c_j]$ if $i\not=j$
and
$$
\sum_{i,j\geq 1} \Big(\mu_j (1-|a_j|)+
2\nu_i(1-|c_i|)\Big)<\infty
$$
is given by
\[
\prod_{i\geq 1} (b_{[c_i]}(p))^{\nu_i}\prod_{j\geq 1}^\star \prod_{k=1}^{\star \mu_j} (b_{a_{jk}}(p))^{\star \mu_j},
\]
where $a_{11}=a_1$ and $a_{jk}\in [a_j]$ are such that $\alpha_{j+1}\not=\bar{\alpha}_{j}$, $j=1,\ldots, n-1$, if $\alpha_j\in\mathbb H\setminus \mathbb R$,  $k=1,2,3,\ldots, \mu_j$.
\end{Tm}
 \begin{proof}
 The Blaschke product converges and defines a slice hyperholomorphic function by Theorem \ref{productblaschke} and its Corollary \ref{Cyproduct}. Let us consider the product:
\begin{equation}\label{prod}
\prod_{i=1}^{\star  \mu_1} (B_{a_{i1}}(p))=B_{a_{11}}(p)\star B_{a_{12}}(p) \star \ldots \star B_{a_{1 \mu_1}}(p).
\end{equation}
As we already observed in the proof of Proposition 5.10 in \cite{acs2} this product admits a zero at the point $a_{11}=a_1$ and it is a zero of multiplicity $1$ if $n_1=1$; if $n_1\geq 2$, the other zeros are $\tilde a_{12}, \ldots, \tilde a_{1 n_1}$ where $\tilde a_{1j}$ belong to the sphere $[a_{1j}]=[a_1]$.
 Thus $\tilde a_{12}, \ldots, \tilde a_{1 n_1}$ all coincide with $a_1$ which is the only zero of the product (\ref{prod}) and it has multiplicity $\mu_1$.
Let us now consider $r\geq 2$ and
\begin{equation}\label{rfactor}
\prod_{j=1}^{\star \mu_r} (B_{a_{rj}}(p))=B_{a_{r1}}(p) \star \ldots \star B_{a_{r n_r}}(p),
\end{equation}
and set
\[
B_{r-1}(p):= \prod_{i\geq 1}^{\star (r-1)} \prod_{k=1}^{\star  \mu_j} (B_{a_{jk}}(p)).
\]
Then from the formula that relates the $\star$-product to the pointwise product (see Proposition 4.3.22 in \cite{MR2752913}) we have that:
\[
B_{r-1}(p)\star B_{a_{r1}}(p)= B_{r-1}(p) B_{a_{r1}}(B_{r-1}(p)^{-1}pB_{r-1}(p))
\]
has a zero  at
 $a_r$ if and only if $B_{a_{r1}}(B_{r-1}(a_r)^{-1}a_rB_{r-1}(a_r))=0$, i.e. if and only if $a_{r1}=B_{r-1}(a_r)^{-1}a_rB_{r-1}(a_r)$. If $n_r=1$ then $a_r$ is a zero of multiplicity $1$ while if $\mu_r\geq 2$, all the other zeros of the product (\ref{rfactor}) belongs to the sphere $[a_r]$ thus the zero $a_r$ has multiplicity $\mu_r$.
 \end{proof}

 We conclude this section by proving that the operator of multiplication by a Blaschke factor is an isometry. In the proof we are in need of the notion of conjugate of a function $f$. Given a slice hyperholomorphic function $f$ consider its restriction to a complex plane $\mathbb C_I$ and write it, as customary, in the form
 $
 f_I(z)=F(z)+G(z)J
 $
 where $J$ is an element in $\mathbb S$ orthogonal to $I$ and $F,G$ are $\mathbb C_I$-valued holomorphic functions.  Define $f^c(p)={\rm ext}(\overline{F(\bar z)}-G(z)J)$ where the extension operator is defined in \eqref{ext}. Note that if $f(p)=\sum_{n\geq 0} p^na_n$ then $f^c(p)\sum_{n\geq 0} p^n\bar a_n$. We have the following:
 \begin{La}
 Let $f\in\mathbf H_2(\mathbb H_+)$. Then $\|f\|_{\mathbf H_2(\mathbb H_+)}=\|f^c\|_{\mathbf H_2(\mathbb H_+)}$.
 \end{La}
 \begin{proof}
 By definition we have
 \[
 \|f\|^2_{\mathbf H_2(\Pi_{+,I})}=\int_{-\infty}^{+\infty} |f_I(Iy)|^2 dy=\int_{-\infty}^{+\infty} (|F(Iy)|^2+|G(Iy)|^2) dy
 \]
and
  \[
  \begin{split}
 \|f^c\|^2_{\mathbf H_2(\Pi_{+,I})}&=\int_{-\infty}^{+\infty}
 |f_I^c(Iy)|^2 dy=\int_{-\infty}^{+\infty} (|\overline{F(-Iy)}|^2+|G(Iy)|^2) dy\\
 &=\int_{-\infty}^{+\infty} (|{F(-Iy)}|^2+|G(Iy)|^2) dy.
 \end{split}
 \]
 Thus $ \|f\|^2_{\mathbf H_2(\Pi_{+,I})}= \|f^c\|^2_{\mathbf H_2(\Pi_{+,I})}$ and taking the supremum for $I\in\mathbb S$ the statement follows.
 \end{proof}
 \begin{Tm}
 Let $b_a$ be a Blaschke factor. The operator
\[
M_{b_a}\,:\, \,f\mapsto b_a\star f
\]
is an isometry from $\mathbf H_2(\mathbb H_+)$ into itself.
\end{Tm}
\begin{proof}
Recall that, by \eqref{blaschkesimpler}, we can write  $b_a(p)=(\tilde p+\bar a)^{-\star}(\tilde p-a)$  for $\tilde p=\lambda^c(p)^{-1}p\lambda(p)$. Let us set $\tilde p=Iy$ where $I\in\mathbb S$. We have
\[
| b_a(Iy)|=|(Iy+\bar a)^{-1}(Iy-a)|=|-(Iy+\bar a)^{-1}\overline{(Iy+a)}|=1.
\]
Similarly, $| b_a^c(Iy)|=1$.  We now observe that for any two functions $f$ and $g$ we have
$(f\star g)^c= g^c\star f^c$. We prove this equality by showing that the two functions $(f\star g)^c$ and $g^c\star f^c$ coincide on a complex plane (so the needed equality follows from the identity principle). Using the notation introduced above, let us write $f_I(z)=F(z)+G(z)J$ and $g_{I}(z)=H(z)+L(z)J$.
We have
\[
(f\star g)_I(z)=f_I(z)\star g_{I}(z)=(F(z)H(z)-G(z)\overline{L(\bar z}))+(F(z)L(z)+G(z)\overline{H(\bar z)})J
\]
so, by definition of $(f\star g)^c$, we have
\[
(f\star g)^c_I(z)=(\overline{F(\bar z)}\, \overline{H(\bar z)}-\overline{G(\bar z})L(z))-(F(z)L(z)+G(z)\overline{H(\bar z)})J
\]
and
\[
\begin{split}
(g^c\star f^c)_I(z)&=(\overline{H(\bar z)}-L(z)J)\star (\overline{F(\bar z)}-G(z)J)\\
&=(\overline{H(\bar z)}\,\overline{F(\bar z)}- L(z) \overline{G(\bar z)})-(\overline{H(\bar z)}G(z)+L(z)F(z))J
\end{split}
\]
the two expressions coincide since the functions $F,G,H,L$ are $\mathbb C_I$-valued and thus they commute.
To compute $\| b_a\star
f\|_{\mathbf{H}_2(\mathbb H_+)}$, where $ f\in\mathbf{H}_2(\mathbb
H_+)$, we  follow an idea used in \cite{abcs-np1} and we compute $\| (b_a\star f)^c\|_{\mathbf{H}_2(\mathbb H_+)}^2$. Note that
$(f^c\star b_a^c)(x+Iy)=0$ where $f^c(x+Iy)=0$, i.e. on a set of isolated points on $\Pi_{+,I}$ while, if $q=f^c(x+Iy)\not=0$,
$(f^c\star b_a^c)(x+Iy)=f^c(x+Iy) b_a^c(q^{-1}(x+Iy)q)$, see \cite[Proposition 4.3.22]{MR2752913}, where $q^{-1}(x+Iy)q={x+I'y}$, see \cite[Proposition 2.22]{ghs}. Thus we have
$(f^c\star b_a^c)(Iy)=f^c(Iy) b_a^c({I'y})$ almost everywhere and
\[
\begin{split}
\| b_a\star f\|_{\mathbf{H}_2(\mathbb H_+)}^2&=\| (b_a\star f)^c\|_{\mathbf{H}_2(\mathbb H_+)}^2\\
&=\sup_{I\in\mathbb S}\int_{-\infty}^{+\infty} |(f^c\star b_a^c)(Iy)|^2dy\\
&=\sup_{I\in\mathbb S}\int_{-\infty}^{+\infty} |f^c(Iy) b_a^c({{I}'y})|^2dy\\
&=\sup_{I\in\mathbb S}\int_{-\infty}^{+\infty} |f^c(Iy)|^2 | b_a^c({I}'y)|^2 dy\\
&=\sup_{I\in\mathbb S}\int_{-\infty}^{+\infty}  |f^c(Iy)|^2dy\\
&=\| f^c\|_{\mathbf{H}_2(\mathbb H_+)}^2.
\end{split}
\]
By the previous lemma, we have $\| f^c\|_{\mathbf{H}_2(\mathbb H_+)}^2=\| f\|_{\mathbf{H}_2(\mathbb H_+)}^2$
and this concludes the proof.
\end{proof}

Blaschke factors will provide a concrete example of the functions
studied in Section \ref{sec6}, see Example \ref{exbl} there.

\section{The Schauder-Tychonoff fixed point theorem}
\setcounter{equation}{0}

In this section we extend the Schauder-Tychonoff
fixed point theorem to the quaternionic setting.
The proof repeats that of the classical
case given in \cite{DS1}, in fact it is readily seen that the arguments hold also in the quaternionic case, but we include it for the reader's convenience. This results is crucial to prove  an invariant subspace theorem for contractions in a Pontryagin spaces.\smallskip
\subsection{The Schauder-Tychonoff fixed point theorem}
In the sequel we will use a consequence of the  Ascoli-Arzel\`a theorem that we state in this corollary.
\begin{La}[Corollary of Ascoli-Arzel\'a theorem]\label{CorAA}
Let $\mathcal{G}_1$ be a compact subset of a topological group $\mathcal{G}$ and let $\mathcal{K}$ be a bounded subset of the space of continuous functions $\mathcal{C}(\mathcal{G}_1)$. Then $\mathcal{K}$ is conditionally compact if and only if for every $\varepsilon>0$ there is a neighborhood
$\mathcal{U}$ of the identity in $\mathcal G$ such that $|f(t)-f(s)|<\varepsilon$ for every $f\in \mathcal K$ and every pair
$s,t\in \mathcal{G}_1$ with $t\in \mathcal{U}s$.
\end{La}
\begin{proof}
It is Corollary 9 p. 267 in \cite{DS1} and its proof can be obtained in the same arguments.
\end{proof}

\begin{Dn}
We say that a quaternionic topological vector space $\mathcal V$ has the fixed point property if
 for every continuous mapping $T:\mathcal{V}\to \mathcal{V}$ there exists $u\in \mathcal{V}$ such that $u=T(u)$.
\end{Dn}
To show our result we need the following Lemmas:
\begin{La}\label{cube}
Let
 $\mathfrak{C}$  be the subspace of $\ell^2(\mathbb{H})$ defined by
 $$
 \mathfrak{C}=\{ \{\xi_n\}\in \ell^2(\mathbb{H}) \ :\ |\xi_n|\leq 1/n,\  \ \ \forall n\in \mathbb{N}  \}.
$$
Then $\mathfrak C$ has the fixed point property.
\end{La}
\begin{proof}
Let $P_n: \ \mathfrak C\to\mathfrak C$ be the map
$$
P_n(\xi_1,\xi_2,\ldots, \xi_n,\xi_{n+1}\ldots )=(\xi_1,\xi_2,\ldots, \xi_n,0,0,\ldots ).
$$
Then $\mathfrak C_n=P_n(\mathfrak C)$ is homeomorphic to the closed sphere in $\mathcal H\cong\mathbb R^{4n}$.
Let now $T: \ \mathfrak C\to\mathfrak C$ be a continuous map. Then  $P_nT:\mathfrak C_n \to \mathfrak C_n$ is continuous.
Brower theorem implies that there is  a fixed point $\zeta_n\in \mathfrak C_n\subseteq \mathfrak C$ and so
$$
|\zeta_n-T(\zeta_n)|\leq \left(\sum_{i=n+1}^\infty \frac{1}{i^2}\right)^{\frac 12}.
$$
Since $\mathfrak C$ is compact,  then $\{\zeta_n\}$ contains a  subsequence converging to a point which is a fixed point of $T$.
\end{proof}
\begin{La}\label{previous}
Let $\mathcal{K}$ be a compact convex subset of a locally convex linear quaternionic space $\mathcal{V}$ and let $T:\mathcal{K}\to \mathcal{K}$ be continuous. If $\mathcal{K}$ contains at least two points, then there exists a proper closed convex subset $\mathcal{K}_1 \subset \mathcal{K}$ such that $T(\mathcal{K}_1)\subseteq \mathcal{K}_1$.
\end{La}
\begin{proof}
It is possible to assume that $\mathcal{K}$ has the $\mathcal{V}^*$ topology.

We will say that a set of continuous linear functionals $F$ is determined by another set $G$, if for every $f\in F$ and $\varepsilon >0$ there exists a neighborhood
$$
\mathcal{N}(0;\gamma,\delta)=\{v\in \mathcal{V} : |g(v)|<\delta, \ \ \ g\in \gamma\},
$$
where $\gamma$ is a finite subset of $G$ with the property that if $u, v\in \mathcal K$ and $u-v\in \mathcal{N}(0;\gamma,\delta)$ then $|f(Tu)-f(Tv)|<\varepsilon$.
It is clear that if $F$ is determined by $G$, then $g(u)=g(v)$ for $g\in G$ implies that $f(Tu)=f(Tv)$ for $f\in F$.
Each continuous linear functional $f$ is determined by some countable set of functional
$G=\{g_m\}_{m\in \mathbb{N}}$.

Thanks to Lemma \ref{CorAA} the scalar function $f(Tu)$ is uniformly continuous on the compact set $\mathcal{K}$.
Hence for every integer $n$ there is a neighborhood $\mathcal{N}(0;\gamma_n,\delta_n)$ of the origin in $\mathcal{V}$,
given by a set of linear continuous functionals $\gamma_n$ and a $\delta_n>0$, such that if $u$, $v\in \mathcal{K}$ and
 $u-v\in \mathcal{N}(0;\gamma_n,\delta_n)$ then $|f(Tu)-f(Tv)|<1/n$. Let $G=\bigcup_{n=1}^\infty\gamma_n$ then $f$ is determined by $G$. It follows that if $F$ is a countable subset of $\mathcal{V}^*$, there exists a countable subset $G_F$
 of $\mathcal{V}^*$ such that each $f\in F$ is determined by $G_F$.
 We claim that each continuous linear functional $f$ can be included in a countable self-determined set $G$ of functionals. In fact, if $f$ is determined by the countable set $G_1$, let each functional in $G_1$ be determined by the countable set $G_2$; then let each functional in $G_2$ be determined by the countable set $G_3$, and so on. We obtain a sequence $\{G_i\}$ and we set
 $G=\{f\}\cup \cup_{i=1}^\infty G_i$.
 Assume now that $\mathcal{K}$ contains two points $u$, $v$, $u\not=v$ and let $f\in \mathcal{V}^*$ be such that $f(u)\not=f(v)$.
 Let $G=\{g_i\}$ be a countable self-determined set of continuous linear functionals containing $f$.
 Since $\mathcal{K}$ is compact, $g_i(\mathcal{K})$ is a bounded set of scalars for every $i$ and since we can multiply $g_i$
 by a suitable constant we may suppose that $g_i(\mathcal{K})\leq 1/i$. In this case the mapping $H: \mathcal{K}\to \ell^2(\mathbb{H})$, defined by
 $$
 H(k):=[g_i(k)]
 $$
 is a continuous mapping of $\mathcal{K}$ onto a compact convex subset $\mathcal{K}_0$ of the subspace $\mathfrak{C}$  of $\ell^2(\mathbb{H})$.
Then $\mathfrak{C}$ contains  trivially at least two points since there are at least two points in $\ell^2(\mathbb C)$, see \cite{DS1}.
 Consider the mapping
 $$
 T_0=HTH^{-1}: \mathcal{K}_0\to \mathcal{K}_0
 $$
 since $G$ is self determined $T_0$ is single-valued. To see that $T_0$ is continuous, let $b_0\in K_0$
 and $\varepsilon \in (0,1)$. Choose $N$ such that $\sum_{i=N+1}^\infty 1/i^2<\varepsilon$.
 Then $G$ is self-determined, there exists a $\delta >0$ and an $m$ such that if $|g_j(u)-g_j(v)|<\delta$,
 $j=1,\ldots, m$ then
 \begin{equation}\label{stareq}
 |g_i(Tu)-g_i(Tv)|<\sqrt{\varepsilon /N},\ \  \ \ i=1,...N.
 \end{equation}
 Thus if $|b-b_0|<\delta$ and $u$ and $v$ are point in $K$ with $b=[g_i(u)]$ and $b_0=[g_i(v)]$ then
 (\ref{stareq}) holds and
 \[
 \begin{split}
 |T_0(b)-T_0(b_0)|^2&=|HTH^{-1}(b)-HTH^{-1}(b_0)|^2
 \\
 &
 \leq \sum_{i=1}^N|g_i(Tu)-g_i(Tv)|^2+2\sum_{i=N+1}^\infty 1/i^2
 \\
 &
 <3 \varepsilon.
 \end{split}
 \]
So $T_0$ is a continuous mapping of $\mathcal{K}_0$ into itself.
From the fixed point property of $\mathfrak C$, see Lemma \ref{cube},  it follows that $T_0$ has  a fixed point $k_0$.
Thus $$TH(k_0)\subseteq H^{-1}T_0(k_0)=H^{-1}(k_0).$$
Setting $\mathcal{K}_1=H^{-1}(k_0)$ we note that $\mathcal{K}_1$  is a proper closed subset of $\mathcal{K}$, and that $T(\mathcal{K}_1)\subseteq \mathcal{K}_1$
The linearity of $H$ implies that $\mathcal{K}_1$ is convex. This concludes the proof.
\end{proof}

\begin{Tm}[Schauder-Tychonoff]
A  compact convex subset of a  locally convex quaternionic  linear space has the fixed point property.
\end{Tm}
\begin{proof}
By the Zorn lemma there exists a minimal convex subset of $\mathcal{K}_1$ of $\mathcal{K}$ with the property that $T\mathcal{K}_1\subseteq \mathcal{K}_1$.
By Lemma \ref{previous} this minimal subset contains only one point.
\end{proof}

\subsection{An invariant subspace theorem}
\setcounter{equation}{0}

As we explained at the beginning of the section, the Schauder-Tychonoff theorem is now used to prove an invariant subspace theorem for contractions in quaternionic Pontryagin spaces. This
theorem is used in the realization theorems to prove the existence of slice hyperholomorphic extensions of certain functions defined
in a neighborhood of a point on the positive axis. In the complex numbers case, this theorem can be found in
\cite[Theorem 1.3.11]{MR92m:47068}. We also refer to \cite[Notes on chapter 1]{MR92m:47068}  for historical notes
on the theorem.

\begin{Tm}
A contraction in a quaternionic Pontryagin space has a unique maximal invariant
negative subspace, and it is one-to-one on it.
\label{tm:inv}
\end{Tm}

\begin{proof} The proof of \cite{MR92m:47068} carries up to the quaternionic setting, and we recall the main lines for the convenience of
the reader. Let  $A$ be a contraction in the Pontryagin space $\mathcal P$. To prove that $A$ has a maximal negative invariant subspace we first recall a well known fact in the theory of linear fractional transformations (see for instance \cite{Dym_CBMS} for more details).
Let  $\mathcal P= \mathcal P_+[\stackrel{\cdot}{+}]\mathcal P_-$ be a fundamental decomposition of $\mathcal P$. Let
\[
A=\begin{pmatrix}A_{11}&A_{12}\\ A_{21}&A_{22}\end{pmatrix}
\]
be the block decomposition of $A$ along $ \mathcal P_+[\stackrel{\cdot}{+}]\mathcal P_-$. Since $A$ is a contraction, and hence a bicontraction (see \cite{acs3}[Theorem 7.2]) we have
\[
A_{21}A_{21}^*-A_{22}A_{22}^*\le -I,
\]
and it follows that $A_{22}^{-1}$ and $A_{22}^{-1}A_{21}$ are strict contractions. Thus the map
\[
L(X)=(A_{11}X+A_{12})(A_{21}X+A_{22})^{-1}
\]
is well defined, and sends in fact the closed unit ball $\mathcal B_1$ of $\mathbf L(\mathcal P_,\mathcal P_+)$ into itself.
The main point in the proof of the theorem is to show that the map $L$ is continuous in the weak operator topology from  $\mathcal B_1$ into itself. Since $\mathcal B_1$ is compact in this topology (and of course convex) the Schauder-Tychonoff theorem implies that $L$ has a unique
fixed point, say $X$. To conclude one notes (see Theorem \cite[1.3.10]{MR92m:47068}) that the space spanned by the elements
\begin{equation}
\label{spaceneq}
f+Xf, \quad f\in\mathcal P_-
\end{equation}
is then negative. It is maximal negative because $X$ cannot have a kernel (any $f$ such that $Xf=0$ will lead to a strictly positive element
of \eqref{spaceneq}).
\end{proof}

\section{The spaces $\mathcal P(S)$}
\setcounter{equation}{0}
\label{sec:HS}
We now introduce the counterparts of the kernels
\eqref{schurgeneral} in the slice hyperholomorphic setting. In the
quaternionic case signature operators are defined as in the
complex case.  Here we consider real signature operators, that is,
which are unitarily equivalent to an operator of the form
\[
\begin{pmatrix}
I_+&0\\0&-I_{-}
\end{pmatrix}.
\]
It is clear that the S-spectrum is concentrated on $\pm 1$, so if $J$ is a signature operator we define
$\nu_-(J)$ as in the complex case. This follows by simple computations, that is $1\pm 2{\rm Re}( s_0) +|s|^2=0$
which give $\pm 1$.\\
In next result we set
$\mathbf L(\mathcal H)\stackrel{\rm def.}{=} \mathbf L(\mathcal
H,\mathcal H)$ where $\mathcal H$ is a two sided quaternionic Hilbert space.
\begin{Dn}
Let $\mathcal H_1$ and $\mathcal H_2$ be two quaternionic
two-sided Hilbert spaces  and let $J_1\in\mathbf L(\mathcal H_1)$
and $J_2\in\mathbf L(\mathcal H_2)$ be two real signature
operators such that $\nu_-(J_1)=\nu_-(J_2)<\infty$. The $\mathbf
L(\mathcal H_1,\mathcal H_2)$-valued function $S$ slice
hypermeromorphic in an axially symmetric s-domain $\Omega$ which
intersects the positive real line belongs to the class $\mathcal
S_\kappa(J_1,J_2)$ if the kernel
\[
K_S(p,q)=J_2k(p,q)-S(p)\star k(p,q)\star_rJ_1S(q)^*
\]
has $\kappa$ negative squares in $\Omega$, where $k(p,q)$ is defined in (\ref{kernel}).
\end{Dn}

We do not mention the dependence of the class on $\Omega$. As we
will see, every element of these classes has a unique meromorphic
extension to $\mathbb H_+$.\smallskip

To reduce the case of arbitrary signature operators (with same number of
negative squares) to the case of the identity, we define the
Potapov-Ginzburg transform in the present setting. We refer to the book
\cite{MR48:904} for the classical case, even though some formulas are also recalled
in \cite{MR2002b:47144}.\smallskip

We begin with a lemma. A proof in the classical case can be found
in \cite[Lemma 4.4.3, p. 164]{adrs}  (the argument there is based on \cite[Lemma 2.1, p. 20]{adps1})  but we repeat the argument for
completeness. First a remark: a matrix $A\in\mathbb H^{m\times
m}$ is not invertible if and only if there exists
$c\not=0\in\mathbb H^{m}$ such that $c^*A=0$. This fact can be
seen for instance from \cite[Theorem 7, p. 202]{MR880123},  where
it is shown that a matrix over a division ring has row rank equal
to the column rank, or \cite[Corollary 1.1.8]{MR1409610}.\smallskip

\begin{La}
Let $T$ be a $\mathbb H^{m\times m}$-valued function
slice hyperholomorphic in an axially symmetric s-domain $\Omega$ which intersect the
positive real line, and such that the kernel
\[
T(p)\star k(p,q)\star_r T(q)^*-k(p,q)I_m
\]
has a finite number of negative squares, say $\kappa$, in $\Omega$.
Then $T$ is invertible in $\Omega$, with the possible exception of a
countable number of spheres.
\label{lemmedupasse}
\end{La}

\begin{proof}
We first show that $T$ is invertible on $\Omega\cap\mathbb R_+$
with the possible exception of a countable number of points. Let
$x_1,\ldots, x_M$ be zeros of $T$. Then, there exist vectors
$c_1,\ldots, c_M$ such that
\[
c_j^*T(x_j)=0,\quad j=1,\ldots, M.
\]
Thus
\[
m_{jk}=c_j^*K(x_j,x_k)c_k=-\frac{c_j^*c_k}{x_j+x_k}.
\]
To conclude we apply \cite[Lemma 2.1, p. 20]{adps1} to the matrix  with block entries $\chi(m_{jk})$ (which is unitarily
equivalent to the matrix $\chi((m_{jk})_)$)) to see that the
$M\times M$ matrix with $jk$ entry $m_{jk}$ is strictly negative,
and so $M\le k$. The result in \cite[Lemma 2.1, p. 20]{adps1} is
proved for the case of complex numbers, but  extends to the
quaternionic case, as is seen by using the map $\chi$ defined in
\eqref{defchi} and Lemma \ref{la:signa}.
\end{proof}

Let now  $S\in\mathcal
S_\kappa(J_1,J_2)$ and let
\begin{equation}
\label{decompS}
S=\begin{pmatrix}S_{11}&S_{12}\\S_{21}&S_{22}\end{pmatrix}
\end{equation}
be the decomposition of $S$ according to fundamental
decompositions of the coefficient spaces. In the
statement of the following theorem, we denote by $I_{2+}$ the identity of
the positive space in the fundamental decomposition of $\mathcal
H_2$.

\begin{Tm}
Let  $S\in\mathcal
S_\kappa(J_1,J_2)$, defined in an axially symmetric s-domain $\Omega$ intersecting the real positive
axis, and with decomposition \eqref{decompS}. Then the function $S_{22}$ is
$\star$-invertible in $\Omega$, with the possible exception of a countable number of
spheres. Let
\begin{equation}
\label{eq:A}
A(p)=\begin{pmatrix} I_{2+}&S_{12}(p)\\
0&S_{22}(p)\end{pmatrix}\quad and\quad \Sigma(p)=\begin{pmatrix}
S_{11}-S_{12}\star S_{22}^{-\star}\star S_{21}&S_{12}\star S_{22}^{-\star}\\
S_{22}^{-\star}\star S_{21}&S_{22}^{-\star}\end{pmatrix}(p).
\end{equation}
Then,
\begin{equation}
\label{perlette}
\begin{split}
J_2k(p,q)-S(p)\star k(p,q)\star_rJ_1S(q)^*
&=\\
&\hspace{-5cm}=A(p)\star\left(k(p,q)-\Sigma(p)\star k(p,q)\star_r\Sigma(q)^*
\right)
\star_r A(q)^*,
\end{split}
\end{equation}
and the kernel
\begin{equation}
\label{perlette2}
k(p,q)-\Sigma(p)\star k(p,q)\star_r \Sigma(q)^*
\end{equation}
has a finite number of negative squares on the domain of
definition of $\Sigma$ in $\Omega$ and hence has a slice
hyperholomorphic extension to the whole of the right half-space,
with the possible exception of a finite number of spheres.
\label{verberie}
\end{Tm}

The function $\Sigma$ is called the Potapov-Ginzburg transform of
$S$, see e.g. \cite[(i), p. 25]{adps1}.

\begin{proof}[Proof of Theorem \ref{verberie}]
To show that $S_{22}$ is $\star$-invertible, we note that
\[
\begin{pmatrix}0&I\end{pmatrix}\left(J_2k(p,q)-S(p)\star k(p,q)\star_rJ_1
S(q)^*\right)\begin{pmatrix}0\\I\end{pmatrix}=
S_{22}(p)\star k(p,q)\star_r S_{22}(q)^*-k(p,q)I_m.
\]
This last kernel has therefore a finite number of negative
squares, and Lemma \ref{lemmedupasse} allows to conclude that
$S_{22}$ is $\star$-invertible, and the definition of the Potapov-Ginzburg transform
makes sense.\smallskip

When $p\in\Omega\cap\mathbb R_+$, the star product is replaced by
the pointwise product and the \eqref{perlette} then follow from
\cite[p. 156]{adrs}. The case of $p\in\Omega$ follows by
slice hyperholomorphic extension. The claim on the number of
negative squares of \eqref{perlette2} follows
\begin{equation}
\label{perlette3}
\begin{split}
k(p,q)-\Sigma(p)\star k(p,q)\star_r\Sigma(q)^*&=\\
&\hspace{-30mm}=
A(p)^{-\star}\star\left(J_2k(p,q)-S(p)\star k(p,q)\star_rJ_1S(q)^*\right)\star_r(A(q)^*)^{-\star_r},
\end{split}
\end{equation}
and from an application of Proposition \ref{pn51}.
\end{proof}

\begin{Dn}
\label{dnks}
Let $S\in\mathcal S_\kappa(J_1,J_2)$. We denote by
$\mathcal P(S)$ the associated reproducing kernel Pontryagin
space of $\mathcal H_2$-valued functions defined in $\Omega$ and
with reproducing kernel $K_S(p,q)$.
\end{Dn}

\section{Realization for elements in $\mathcal
S_\kappa(J_1,J_2)$}
\setcounter{equation}{0}
\label{sec6}

In this section we present a realization theorem for elements in
$\mathcal S_\kappa(J_1,J_2)$, where the state space is the
reproducing kernel Pontryagin space $\mathcal P(S)$ (see
Definition \ref{dnks} for the latter). In the case $\kappa=0$ one
could get the existence of a realization using a Cayley transform
in the variable and use our previous results in \cite{acs1}.
Here we give a direct proof to get a realization defined in
$\mathcal P(S)$, taking into account that $\kappa$ may be strictly
positive. We begin with a definition:

\begin{Dn}
Let $\mathcal P_1$ and $\mathcal P_2$ be two quaternionic right Pontryagin spaces.
A pair of operators $(G,A)\in\mathbf L(\mathcal P_1,\mathcal P_2)\times\mathbf L(\mathcal P_1)$
is called observable (or closely outer connected) if
\[
\cap_{n=0}^\infty\ker G A^n=\left\{0\right\}.
\]
\end{Dn}

The terminology {\sl observable} is the one from the theory of linear systems, while  {\sl  closely outer connected} has been used in
operator theory in particular by Krein and Langer, see \cite{adrs}.

\begin{Tm}
Let $x_0$ be a strictly positive real number. A function $S$
slice hyperholomorphic in an axially symmetric s-domain $\Omega$ containing $x_0$ is
the restriction to $\Omega$ of an element of $\mathcal
S_\kappa(J_1,J_2)$ if and only if it can be written as
\begin{equation}
\label{realS}
\begin{split}
S(p)&=H-(p-x_0)\left(G-(\overline{p}-x_0)(\overline{p}+x_0)^{-1}GA\right)\times \\
&\hspace{10mm}\times
\left(\frac{|p-x_0|^2}{|p+x_0|^2}A^2 - 2{\rm
Re}~\left(\frac{p-x_0}{p+x_0}\right)A+I\right)^{-1}F,
\end{split}
\end{equation}
where $A$ is a linear bounded operator in a right-sided
quaternionic Pontryagin space $\Pi_\kappa$ of index $\kappa$,
and, with $B=-(I+x_0A)$, the operator matrix
\[
\begin{pmatrix}
B&F\\ G&H\end{pmatrix}\,\,:\,\, \begin{pmatrix}\Pi_k\\
\mathcal H_1\end{pmatrix}\,\,\,\longrightarrow
\,\,\,\, \begin{pmatrix}\Pi_k\\ \mathcal H_2\end{pmatrix}
\]
is co-isometric. In particular $S$ has a unique slice
hypermeromorphic extension to $\mathbb H_+$. Furthermore, when
the pair $(G,A)$ is observable, the realization is unique up to a
unitary isomorphism of Pontryagin right quaternionic spaces.
\label{thmrealS}
\end{Tm}

\begin{Rk} {\rm When the operators are finite matrices we note that formula
\eqref{realS} can be rewritten as:
\[
S(p)=H-(p-x_0) G\star((x_0+p)I+(p-x_0)B)^{-\star} F.
\]
Sometimes, and by abuse of notation, we will use this expression
also for the infinite dimensional case, see Proposition
\ref{formula060813} for more information.}
\end{Rk}

\begin{proof}[Proof of Theorem \ref{thmrealS}]
We proceed in a number of steps, and first prove in Steps 1-8
that a realization of the asserted type exists with
$\Pi_k=\mathcal P(S)$. We denote by $ \mathcal H_2({J_2})$ the
space $ \mathcal H_2$ endowed with the indefinite inner product
\[
[u,v]_{J_2}=[u,J_2v]
\]
and similarly $ \mathcal H_1({J_1})$. Both $ \mathcal H_1({J_1})$
and $ \mathcal H_2({J_2})$  are quaternionic
Pontryagin spaces, and they have the same index.\\

Following \cite[pp. 51-52]{abds2-jfa} we introduce a relation $R$
in $(\mathcal P(S)\oplus \mathcal H_2({J_2}))\times(\mathcal
P(S)\oplus\mathcal H_1( {J_1}))$ by the linear span of the vectors
\begin{equation}
\left(
\begin{pmatrix}K_S(\cdot, q)(x_0-\overline{q})u\\
(x_0-\overline{q})v\end{pmatrix}\,,\,
\begin{pmatrix}
K_S(\cdot, q)(x_0+\overline{q})u-2x_0K_S(\cdot,
x_0)u+\sqrt{2x_0}K_S(\cdot, x_0)(x_0-\overline{q})v\\
\sqrt{2x_0}(S(q)^*-S(x_0)^*)u+
S(x_0)^*(x_0-\overline{q})v\end{pmatrix} \right).
\end{equation}

STEP 1: {\sl The relation $R$ is isometric.}\\

Indeed, let $(F_1,G_1)$ and $(F_2,G_2)$ be two elements in the
relation, corresponding to $q_1\in\Omega, u_1,v_1\in\mathcal H_1$
and to $q_2\in\Omega, u_2,v_2\in\mathcal H_2$ respectively. On
the one hand we have
\[
\begin{split}
[F_2,F_1]&=[(x_0-q_1)K_S(q_1,q_2)(x_0-\overline{q_2})u_2,u_1]+
[(x_0-q_1)J_2(x_0-\overline{q_2})v_2,v_1].
\end{split}
\]
On the other hand, with $G_1=\begin{pmatrix}g_1\\
h_1\end{pmatrix}$ where
\[
\begin{split}
g_1(\cdot)&=K_S(\cdot,
q_1)(x_0+\overline{q_1})u_1-{2x_0}K_S(\cdot,
x_0)u_1+\sqrt{2x_0}K_S(\cdot, x_0)(x_0-\overline{q_1})v_1,\\
h_1&=
\sqrt{2x_0}(S(q_1)^*-S(x_0)^*)u_1+S(x_0)^*(x_0-\overline{q_1})v_1
\end{split}
\]
(and similarly for $G_2$) we have
\[
\begin{split}
[G_2,G_1]=[g_2,g_1]+[h_2,h_1].
\end{split}
\]
We want to show that
\begin{equation}
\label{FFGG}
[F_2,F_1]=[g_2,g_1]+[h_2,h_1].
\end{equation}
In the computations of these inner products, there are terms
which involve only $u_1,u_2$, terms which involve only $v_1,v_2$
and similarly for $u_1,v_2$ and $v_1,u_2$. We now write these
inner terms separately:\\

\underline{Terms involving $u_1,u_2$.} To show that these terms
are the same on both sides of \eqref{FFGG} we have to check that
\[
\begin{split}
[({x_0}-q_1)K_S(q_1,q_2)({x_0}-\overline{q_2})u_2,u_1]&=
[({x_0}+q_1)\left(K_S(q_1,q_2)({x_0}+\overline{q_2})-\right.\\
&\hspace{5mm} -2x_0K_S({x_0},q_2)({x_0}+
\overline{q_2})-\\
&\hspace{5mm}\left.-2x_0({x_0}+q_1)K_S(q_1,{x_0})
+4x_0^2K_S({x_0},{x_0})\right)u_2,u_1]+\\
&\hspace{5mm}+2x_0[(S(q_1)-S({x_0}))J_1(S(q_2)^*-S({x_0}))u_2,u_1].
\end{split}
\]
Using
\begin{equation}
\label{eq:k11} k(x_0,x_0)=\frac{1}{2x_0}\quad{\rm and}\quad
K_S(x_0,x_0)=\frac{1}{2x_0}\left(J_2-S(x_0)J_1S(x_0)^*\right),
\end{equation}
we see that this is equivalent to prove that
\[
\begin{split}
(x_0-q_1)J_2k(q_1,q_2)(x_0-\overline{q_2})-
(x_0-q_1)S(q_1)J_1k(q_1,q_2)S(q_2)^*(x_0-\overline{q_2})&=\\
&\hspace{-11cm}= (x_0+q_1)J_2k(q_1,q_2)(x_0+\overline{q_2})-
(x_0+q_1)S(q_1)J_1k(q_1,q_2)S(q_2)^*(x_0+\overline{q_2})-\\
&\hspace{-10cm}-2x_0(J_2-S(q_1)J_1S(x_0)^*)-2x_0(J_2-S(x_0)J_1S(q_2)^*)+\\
&\hspace{-10cm}+
2x_0\left(J_2-S(x_0)J_1S(x_0)^*\right)+\\
&\hspace{-10cm}+2x_0(S(q_1)-S(x_0)J_1(S(q_2)^*-S(x_0)^*).
\end{split}
\]
But this amounts to check that
\[
q_1k(q_1,q_2)+k(q_1,q_2)\overline{q_2}=1,
\]
which has been seen to hold in Proposition \ref{identity}.\\

\underline{Terms involving $v_1,v_2$.} To show that these terms
are the same on both sides of \eqref{FFGG} we have to check that
\[
\begin{split}
(x_0-q_1)J_2(x_0-\overline{q_2})&=
[(x_0-q_1)S(x_0)J_1S(x_0)^*(x_0-\overline{q_2})v_2,v_1]+\\
&\hspace{5mm}+2x_0[(x_0-q_1)K_S(x_0,x_0)(x_0-\overline{q_2})v_2,v_1].
\end{split}
\]
This follows directly from the formula for $K_S(x_0,x_0)$, see
\eqref{eq:k11}.\\

\underline{Terms involving $u_2,v_1$.} There are no such terms on
the left side of \eqref{FFGG} and so we need to show that the
terms on the right add up to $0$. This is the case since
\[
\begin{split}
&\hspace{5mm}\sqrt{2x_0}[(x_0-q_1)S(x_0)J_1(S(q_2)^*-S(x_0)^*)u_2,v_1]+\\
&\hspace{5mm}+\sqrt{2x_0}[(x_0-q_1)\left(
K_S(x_0,q_2)(x_0+\overline{q_2})-
2x_0K_S(x_0,x_0)\right)u_2,v_1]=\\
&=[Xu_2,v_1]\\
&=0
\end{split}
\]
with
\[
\begin{split}
X&=\sqrt{2x_0}(x_0-q_1)S(x_0)J_1(S(q_2)^*-S(x_0)^*)+\sqrt{2x_0}(x_0-q_1)
\left(J_2-S(x_0)J_1S(q_2)^*\right)-\\
&\hspace{5mm} -\sqrt{2x_0}(x_0-q_1)\left(J_2
-S(x_0)J_1S(x_0)^*\right)\\
&=0
\end{split}
\]
since
\[
K_S(x_0,q_2)(x_0+\overline{q_2})=J_2-S(x_0)J_1S(q_2)^*.
\]

\underline{Terms involving $u_1,v_2$.} These form a symmetric
expression to the previous one, and will not be written down.\\

STEP 2: {\sl The domain of $R$ is dense.}\\

To prove this step, let $\begin{pmatrix}f\\ w\end{pmatrix}\in
(\mathcal P(S)\oplus \mathcal H_2({J_2}))$ be orthogonal to ${\rm
Dom}\,R$. Then, for all $q\in\Omega$ and $u,v\in\mathcal H_2$ we
have
\[
[(x_0-q)f(q),u]+[(x_0-q)w,v]_{J_2}=0.
\]
It follows that $w=0$ and that
\[
(x_0-q)f(q)\equiv 0,\quad q\in\Omega,
\]
and so $f\equiv0$ in $\Omega$.\\

STEP 3: {\sl The relation $R$ extends to the graph of an isometry}\\

Indeed, the spaces  $\mathcal P(S)\oplus \mathcal H_2({J_2})$ and
$\mathcal P(S)\oplus\mathcal H_1( {J_1})$ are Pontryagin spaces
with same index. By the quaternionic version of a theorem of
Shumlyan (see \cite[Theorem 1.4.1, p. 27]{adrs} for the classical
case and \cite[Theorem 7.2]{acs2} for the quaternionic case) a
densely defined contractive relation defined on a pair of
Pontryagin spaces with same index extends to the
graph of a contraction.\\

In preparation to the next step we introduce an operator
$R_{x_0}$ as follows. Let $\mathcal H$ be a two-sided
quaternionic Hilbert space. A $\mathcal H$-valued function slice
hyperholomorphic in a neighborhood of $x_0>0$ can be written as
a convergent power series
\[
f(p)=\sum_{n=0}^\infty (p-x_0)^nf_n,
\]
where the coefficients $f_n\in\mathcal H$. We define
\begin{equation}
(R_{x_0}f)(p)=(p-x_0)^{-1}(f(p)-f(x_0))\stackrel{\rm def.}{=}
\begin{cases}\,\,\sum_{n=1}^\infty (p-x_0)^{n-1}f_n,\quad
p\not=x_0,\\
\quad f_1,\hspace{2.8cm}\,\,\, p=x_0.
\end{cases}
\label{RX0}
\end{equation}

STEP 4: {\sl Let $V$ denotes the isometry in the previous step.
We compute $V^*$ and show that, with
\begin{equation}
\label{V*}
V^*=\begin{pmatrix}B&F\\ G& H\end{pmatrix}\,\,:\,\,\,\,\, \mathcal
P(S)\oplus \mathcal H_2({J_2})\,\,\Longrightarrow\,\, \mathcal
P(S)\oplus \mathcal H_2({J_1}),
\end{equation}
we have $H=S(x_0)$ and

\begin{align}
\label{I1}
Bf&=-(I+2x_0R_{x_0})f,\\
\label{I2}
Fu&=-\sqrt{2x_0}R_{x_0}Su,\\
\label{I3}
Gf&=\sqrt{2x_0}f(x_0).
\end{align}
}

To compute \eqref{I1} let $f\in\mathcal P(S)$ and
$(p,u)\in\Omega\times\mathcal H_2$. We have
\[
\begin{split}
[(x_0-p)(Bf(p)),u]&= [Bf,K_S(\cdot,
p)(x_0-\overline{p})u]\\
&=[f,B^*(
K_S(\cdot, p)(x_0-\overline{p})u)]\\
&=[f,K_S(\cdot,
p)(x_0+\overline{p})u-2x_0K_S(\cdot, x_0)u]\\
&=[(p+x_0)f(p)-2x_0f(x_0),u],
\end{split}
\]
and so
\[
(x_0-p)(Bf(p))=(p+x_0)f(p)-2x_0f(x_0),\quad p\in\Omega,
\]
which can be rewritten as \eqref{I1}.\\

 Similarly, to compute \eqref{I2} let
$v\in\mathcal H_2$. We have:
\[
\begin{split}
[(x_0-p)((Fv)(p)),u]&= [Fv,K_S(\cdot,
p)(x_0-\overline{p})u]\\
&=[v,\sqrt{2x_0}(S(p)^*-S(x_0)^*)u]\\
&=[\sqrt{2x_0}(S(p)-S(x_0))v,u],
\end{split}
\]
and so
\[
(x_0-p)(Fv(p))=\sqrt{2x_0}(S(p)-S(x_0))v,\quad p\in\Omega.
\]
Finally, we have:
\[
\begin{split}
[(x_0-p)Gf,v]&=
[Gf,(x_0-\overline{p})v]\\
&=[f,G^*(x_0-\overline{p})v)]\\
&=[f,\sqrt{2x_0}K_S(\cdot, x_0)(x_0-\overline{p})v]\\
&=\sqrt{2x_0}[(x_0-p)f(x_0),v].
\end{split}
\]
where we have used \eqref{rightcondition} to get the first
equality.\\

STEP 5: {\sl We prove \eqref{realS} for real $p$ near $x_0$:} The operator $I+2x_0R_{x_0}$ is bounded and so is the
operator $R_{x_0}$ (with $x_0>0$).
Let $f\in\mathcal P(S)$, with power series expansion
\[
f(p)=\sum_{n=0}^\infty (p-x_0)^nf_n,\quad f_0,f_1,\ldots \in\mathcal H_2,
\]
around $x_0$. We have for real $p=x$ near $x_0$:
\[
\begin{split}
f(x)&=\sum_{n=0}^\infty (x-x_0)^nf_n\\
&=\frac{1}{\sqrt{2x_0}}\sum_{n=0}^\infty (x-x_0)^nGR_{x_0}^nf\\
&=\frac{1}{\sqrt{2x_0}}G(I-(x-x_0)R_{x_0})^{-1}f.
\end{split}
\]
Applying this formula to $f=R_{x_0}Su=-\frac{1}{\sqrt{2x_0}}Fu$ where $u\in\mathcal H_1$ we have
\[
\begin{split}
(R_{x_0}Su)(x)&=-G(2x_0I-2(x-x_0)x_0R_{x_0})^{-1}Fu
\end{split}
\]
and so, since $B=-I-2x_0R_{x_0}$,
\[
\begin{split}
S(x)u&=S(x_0)u+(x-x_0)(R_{x_0}Su)(x)\\
&=S(x_0)u-(x-x_0) G(2x_0I-2(x-x_0)x_0R_{x_0})^{-1}Fu\\
&=S(x_0)u-(x-x_0) G(2x_0I+(x-x_0)(B+I))^{-1}Fu\\
&=S(x_0)u-(x-x_0) G((x+x_0)I+(x-x_0)B)^{-1}Fu.
\end{split}
\]

STEP 6: {\sl Assume that $J_1=I_{\mathcal H_1}$ and $J_2=I_{\mathcal H_2}$.
Then, the operator $(x_0+x)I+(x-x_0)B$ is invertible for all real $x$, with the possible exception of a finite set in $\mathbb R$.}\\

Assume first the kernel $K_S$ to be positive definite. Then, the operator matrix \eqref{V*} is a contraction
between Hilbert spaces and so $B$ is a Hilbert space contraction, and the operator
\[
I-\frac{x_0-x}{x_0+x}B
\]
is invertible for all $x>0$, with the possible exception of a finite set,  since $|\frac{x_0-x}{x_0+x}|<1$ for such $x$.\smallskip

Assume now that $\mathcal P(S)$ is a Pontryagin space. The operator $V^*$ is a contraction between Pontryagin
spaces of same index, and so its adjoint $V$ is a contraction (see  \cite[Theorem 7.2]{acs3}).
So it holds that
\[
B^*B+G^*G\le I.
\]
But
\[
\langle G^*Gf,f\rangle=\langle Gf,Gf\rangle_{\mathcal H_2}\ge 0
\]
since $J_2=I_{\mathcal H_2}$ and so $B$ is a contraction. It
admits a maximal strictly negative invariant subspace, say
$\mathcal M$ (see \cite[Theorem 1.3.11]{MR92m:47068} for the
complex case and Theorem \ref{tm:inv} for the quaternionic case).
Writing
\[
\mathcal P(S)=\mathcal M[+]\mathcal M^{[\perp]},
\]
the operator matrix representation of $B$ is upper triangular with
respect to this decomposition where
\[
B=\begin{pmatrix} B_{11}&B_{12}\\0&B_{22}\end{pmatrix}.
\]
The operator $B_{22}$ is a contraction from the
Hilbert space $\mathcal M^{[\perp]}$ into itself, and so
$I-\frac{x_0-x}{x_0+x}B_{22}$ is invertible for every $x>0$, with the possible exception of a finite set.
The operator $B_{11}$ is a contraction from the finite dimensional anti-Hilbert
space $\mathcal M$ onto itself, and so has right eigenvalues {\sl outside
the open unit ball}. So the operator $I-\frac{x_0-x}{x_0+x}B_{11}$,
is invertible in $x>0$, except the points $x\not=x_0$
such that $\frac{x+x_0}{x-x_0}$ is a real eigenvalue of $B_{11}$
of modulus greater or equal to $1$. There is a finite number of such points
since, see \cite[Corollary 5.2, p. 39]{MR97h:15020}, a $n\times n$ quaternionic matrix
has exactly $n$ right eigenvalues (counting multiplicity)
up to equivalence (in other words, it has exactly $n$ spheres of eigenvalues).\smallskip

It follows that the operator
\[
I-\frac{x_0-x}{x_0+x}B=\begin{pmatrix}
I-\frac{x_0-x}{x_0+x}B_{11}
&-\frac{x_0-x}{x_0+x}B_{12}\\0&
I-\frac{x_0-x}{x_0+x}
B_{22}\end{pmatrix}
\]
is invertible for all $x>0$, with the possible exception of a finite number of points.\\

STEP 7: {\sl Assume that $J_1=I_{\mathcal H_1}$ and $J_2=I_{\mathcal H_2}$.
The function $S$ admits a slice hypermeromorphic
extension to $\mathbb H_+$, with the possible exception of a finite number of spheres.}\\

We note that, for $p\in\mathbb H$ near $x_0$  we can extend $S(x)u$ computed in STEP 5 to a slice hyperholomorphic function:
\[
\begin{split}
S(p)u&=S(x_0)u+\frac{x_0-p}{x_0+p} G \star \left(I-\frac{x_0-p}{x_0+p}B\right)^{-\star}Fu\\
=&S(x_0)u+\frac{p-x_0}{p+x_0}\star\left(G- \frac{x_0-\bar p}{x_0+\bar p}
GB\right)\left(\frac{|x_0-p|^2}{|x_0+p|^2}B^2 - 2{\rm
Re}~\left(\frac{x_0-p}{x_0+p}\right)B+I\right)^{-1}Fu.
\end{split}
\]

Let $t=\frac{{\rm Re}\,q}{|q|^2}$ where $q=\frac{x_0-p}{x_0+p}$. We have

\[
\begin{split}
|q|^2B^2-2({\rm Re}\,q)B+I&=
|q|^2\begin{pmatrix}B_{11}^2-2tB_{11}+\frac{1}{|q|^2}&
B_{11}B_{12}+B_{12}B_{22}-2tB_{12}+\frac{1}{|q|^2}\\
0&B_{22}^2-2tB_{22}+\frac{1}{|q|^2}\end{pmatrix}.
\end{split}
\]

By the property of the resolvent, the operator
$B_{22}^2-2tB_{22}+\frac{1}{|q|^2}$ is invertible for $q$ such
that $\frac{1}{|q|^2}$ is in the resolvent set of $B_{22}$. Since
$B_{22}$ is a Hilbert space contraction, this happens in
particular when $|q|<1$, see \cite{MR2752913}, proof of Theorem 4.8.11. Similarly  the operator
$B_{11}^2-2tB_{11}+\frac{1}{|q|^2}$ is invertible if and only if
$\frac{1}{|q|^2}$ is in the resolvent set of $B_{11}$. Since $B_{11}$
is a finite dimensional Hilbert space expansion, it has just point S-spectrum which is
inside the closed unit ball. The point S-spectrum coincides with the set of right eigenvalues, see Remark \ref{spettro1}, and  it consists of a finite number of (possibly degenerate) spheres.\\

We now consider the case of arbitrary signature
matrices, with same negative index.\\

STEP 8: {\sl We use the Potapov-Ginzburg transform to show that
$S$ has a
meromorphic extension.}\\

This follows from computing $S$ from its Potapov-Ginzburg
transform.\\

STEP 9: {\sl Any $S$ with a realization of the form \eqref{realS}
is in a class  $\mathcal S_\kappa(J_1,J_2)$.}\\

Indeed, for real $p=x$ and $q=y$ near $x_0$, the existence of the
realization leads to
\[
\frac{J_2-S(x)J_1S(y)^*}{x+y}=G(I(x_0+x)-(x+x_0)B)^{-1}(I(y+x_0)-(y-x_0)B)^{-*}G^*,
\]
where $B=-(I+x_0A)$. Thus, with
$K(x,y)=G(I(x_0+x)-(x+x_0)B)^{-1}(I(y+x_0)-(y-x_0)B)^{-*}G^*$,
\[
J_2-S(x)J_1S(y)^*=xK(x,y)+K(x,y)y
\]
and the result follows by observing that \eqref{realS} is the hyperholomorphic extension.\\

STEP 10: {\sl The realization is unique up to isomorphism when it
is observable.}\\

We follow \cite{adrs}. Let $p$ be a real number and set
$x=\frac{p-x_0}{p+x_0}$. When $p$ varies in a real neighborhood
of $x_0$ then $x$ varies in a real neighborhood $\mathcal{I}_0$ of the origin.
For $x,y\in \mathcal{I}_0$ we have \[
\frac{J_2-S(x)J_1S(y)^*}{1-xy}=G_1(I_{\mathcal
P_1}-xB_1)^{-1}(I_{\mathcal P_1}-yB_1)^{-*}G_1^*= G_2(I_{\mathcal
P_2}-xB_2)^{-1}(I_{\mathcal P_2}-yB_2)^{-*}G_2^*,
\]
where the indices $1$ and $2$ correspond to two observable and
coisometric realizations, with state spaces $\mathcal P_1$ and
$\mathcal P_2$ respectively. Then the domain and range of the
relation $R$ spanned by the pairs
\[
((I_{\mathcal P_1}-yB_1)^{-*}G_1^*h,(I_{\mathcal P_2}-yB_2)^{-*}G_2^*k),
\quad h,k\in\mathcal H_2,
\]
are dense. By the quaternionic version of a theorem of Shmulyan
(see \cite[Theorem 7.2]{acs2}) $R$ is the graph of a unitary map,
which provides the desired equivalence. The arguments are as in
\cite{adrs}.
\end{proof}
\begin{Rk}{\rm
In the case $x_0\rightarrow 0$, we can rewrite the computations with
\[
\begin{split}
G_0f&=f(x_0)\\
F_0u&=R_{x_0}Su\\
I+B&=-x_0R_{x_0},\\
I-B&=2I+x_0R_{x_0};
\end{split}
\]
we have that
\[
\begin{split}
S(p)u&=S(x_0)u+(p-x_0) G \star\left(\mathcal
M^\ell_{x_0+{p}}+\mathcal M^\ell_{p-x_0}B\right)^{-1}F_0u\\
&=S(x_0)+(p-x_0)2x_0
G_0\star\left(x_0(I-B)+p(I+B)\right)^{-1}F_0u\\
&=S(x_0)+(p-x_0) G_0\star\left(I+x_0R_{x_0}-pR_{x_0}\right)^{-1}F_0u\\
\end{split}
\]
which tends formally to the backward shift realization as
$x_0\rightarrow 0$.}
\label{remarkbw}
\end{Rk}
\begin{Ex}{\rm We now show how to obtain a realization for a Blaschke factor $b_a(p)$.
For real $p=x$, using formula \eqref{blaschkesimpler} we obtain that $b_a(x)=(x-\bar a)^{-1}(x-a)$, moreover
\[
\begin{split}
b_{a}(x)&=b_{a}(1)+b_{a}(x)-b_{a}(1)\\
&=\frac{1-a}{1+\bar a}+(x-1)\frac{2{\rm
Re}(a)}{(x+\overline{a})(1+\overline{a})}\\
&= \frac{1-a}{1+\bar a}+(x-1)\frac{2{\rm
Re}(a)}{(x+\frac{1-B}{1+B})(1+\overline{a})},\quad{\rm
where}\,\,
B=\frac{1-\overline{a}}{1+\overline{a}}\\
&=\frac{1-a}{1+\bar a}+(x-1)\frac{2 {\rm
Re}(a)(1+B)}{(x(1+B)+(1-B))(1+\overline{a})}\\
&=\frac{1-a}{1+\bar a}+(x-1)\frac{2{\rm
Re}(a)}{(x(1+B)+(1-B))}\frac{2}{(1+\overline{a})^2}
\\
&=\frac{1-a}{1+\bar a}+(x-1)\frac{2{\rm
Re}(a)}{1+\overline{a}}((x+1)+(x-1)B)^{-1}\frac{2}{1+\overline{a}}
\end{split}
\]
since
\[
\frac{1+B}{1+\overline{a}}=\frac{2}{(1+\overline{a})^2}.
\]
Now note that
\[
b_{a}(p)=H-(p-1) G\star ((p+1)+(p-1)B)^{-\star} F
\]
is slice  hyperholomorphic, extends $S(x)$, and
\[
\begin{pmatrix}B&F\\G&H\end{pmatrix}=\begin{pmatrix}
\dfrac{1-\overline{a}}{1+\overline{a}}&\dfrac{2\sqrt{{\rm Re}\,
a}}{1+\overline{a}}\\
\\
-\dfrac{2\sqrt{{\rm Re}\, a}}{1+\overline{a}}& \dfrac{1-a}{1+\bar a}
\end{pmatrix}.
\]
This matrix is unitary.

}
\label{exbl}
\end{Ex}

We now present an example of functions in a class $\mathcal S_0(J,J)$. Consider a linear bounded operator $A$ in a right
quaternionic Hilbert space $\mathcal H$, and assume that $A+A^*$ is finite dimensional, say of rank $m$. We can thus write:
\[
A+A^*=-CJC^*,
\]
where $J\in\mathbb H^{m\times m}$ is a real signature matrix, and where $C$ is linear
bounded operator from $\mathbb H^m$ into
$\mathcal H$. We will assume $(C,A)$ observable.
\begin{Rk}
The pair $(C,A)$ is observable if and only if there is no non
trivial invariant subspace of $A$ on which $A+A^*=0$.
\end{Rk}

The proof of this lemma is as in the complex case, and will be
omitted.\\

We conclude this section with an example of a function in $\mathcal S_0(J,J)$, which, by analogy with the classical case, we call the
characteristic operator function of the operator. Connections with operator models will be considered elsewhere, but we remark here that the function $S$ in
\eqref{COF} defined uniquely $A$ when the pair $(C,A)$ is observable.

\begin{Dn}
The function
\begin{equation}
\label{COF}
S(p)=I-p C^*\star(I-pA)^{-\star}CJ
\end{equation}
is called the characteristic operator function of the operator $A$.
\end{Dn}

\begin{Tm}
The characteristic operator function belongs to $\mathcal S_0(J,J)$.
\end{Tm}

\begin{proof}
Let
\[
K(p,q)=C^*\star
(I-pA)^{-\star}(I-\overline{q}A^*)^{-\star_r}\star_rC.
\]
Then it holds that
\[
J-S(p)JS(q)^*=p K(p,q)+K(p,q)\overline{q}.
\]
This formula is proved by first considering the case of real $p$
and $q$, and taking the slice hyperholomorphic extension, and
proves that $S\in\mathcal S_0(J,J)$.
\end{proof}

We note that formula \eqref{COF} corresponds to a realization
centered at $0$, as in Remark \ref{remarkbw}, and not to a
realization of the form \eqref{realS}. It would be interesting to find
a functional model for the operator $A$ in terms of $S$. The
special case where $S$ is a (possibly infinite convergent)
Blaschke product is of special interest. The case of general $S$
leads to the question of finding the $\star$-multiplicative
structure of elements in $\mathcal S_0(J,J)$, that is the
counterpart of the paper \cite{pootapov} in the present setting.

\section{The space $\mathcal L(\Phi)$ and realizations for generalized
positive functions}
\label{Sec:2}

In the present section we give realization for a generalized positive function with  $\mathcal L(\Phi)$ as state space. Note that
a Cayley transform (with real coefficients) will map a generalized positive function into a generalized Schur function, and even more a Cayley trasform on the variable will reduce the problem to the case of a Schur function of the quaternionic unit ball. But this procedure will not lead an intrinsic
realization in the natural space associated to generalized positive function.

\setcounter{equation}{0}
\subsection{The indefinite case}
\begin{Dn}
Let $\mathcal H$ be a quaternionic Hilbert space, and let
$J\in\mathbf L(\mathcal H)$ be a real signature operator. A $\mathbf
L(\mathcal H)$-valued function $\Phi$ slice hyperholomorphic in
an axially symmetric s-domain $\Omega$ which intersects the
positive real line belongs to the class ${\rm GP}_\kappa(J)$ if
the kernel
\begin{equation}
\label{rkpsphi}
K_\Phi(p,q)=J\Phi(p)\star k(p,q)+
k(p,q)\star_r\Phi(q)^*J
\end{equation}
has $\kappa$ negative squares in $\Omega$.
\end{Dn}

\begin{La}
The kernel $K_\Phi$ satisfies
\begin{equation}
\label{equaphi}
pK_\Phi(p,q)+K_\Phi(p,q)\overline{q}=J\Phi(p)+\Phi(q)^*J.
\end{equation}
\end{La}
\begin{proof}
It follows with immediate computations from Proposition \ref{identity}.
\end{proof}

As in the case of generalized Schur functions, we do not mention
the dependence of the class on $\Omega$ since, as we prove later,
every element of a class ${\rm GP}_\kappa(J)$ has a unique
slice hypermeromorphic extension to $\mathbb H_+$.\\

We note that $J$ does not play a role, as noted in \cite[p. 358,
footnote]{kl1}, and could be set to be the identity. We denote by
$\mathcal H$ a two sided quaternionic Hilbert space, and recall that
$\mathbf L(\mathcal H){=} \mathbf L(\mathcal
H,\mathcal H)$.

\begin{Tm}
A $\mathbf L(\mathcal H)$-valued function $\Phi$ slice
hyperholomorphic in an axially symmetric s-domain $\Omega$
containing $x_0>0$ is in the class ${\rm GP}_\kappa(J)$ if and
only if there exists a right quaternionic Pontryagin space
$\Pi_\kappa$ of index $\kappa$ and operators
\[
\begin{pmatrix}B&F\\ G&H\end{pmatrix}\,\,:\,\, \begin{pmatrix}\Pi_k\\
\mathcal H\end{pmatrix}\,\,\,\longrightarrow
\,\,\,\, \begin{pmatrix}\Pi_k\\ \mathcal H\end{pmatrix}
\]
verifying
\[
\begin{split}
(I+2x_0B)(I+2x_0B)^*&=I\\
\end{split}
\]
and such that $\Phi$ can be written as
\begin{equation}
\Phi(p)=H-(p-x_0) G\star((p+x_0)I+(p-x_0)B)^{-\star} F.
\label{eq:realphi12345}
\end{equation}
Furthermore, $\Phi$ has a unique slice hypermeromorphic extension
to $\mathbb H_+$.  Finally, when the pair $(G,B)$ is observable,
the realization is unique up to a unitary isomorphism of
Pontryagin right quaternionic spaces. \label{thmrealP}
\end{Tm}

\begin{proof} Given $\Phi\in{\rm GP}_\kappa(J)$, we denote by $\mathcal L(\Phi)$
associated right reproducing kernel Pontryagin space of $\mathcal
H$-valued functions with reproducing kernel $K_\Phi$. We
proceed in a number of steps to prove the theorem.\\

STEP 1: {\sl The formula
\begin{equation}
(p-x_0)(Bh(p))=(p+x_0)h(p)-2x_0h(x_0),\quad h\in\mathcal L(\Phi).
\end{equation}
defines a (continuous) coisometry in $\mathcal L(\Phi)$.}\\

Indeed, define a relation $\mathcal R_{x_0}$ on $\mathcal L(\Phi)\times\mathcal
L(\Phi)$ generated by the linear span of the pairs
\begin{equation}
\mathcal R_{x_0}=\left(K_\Phi(\cdot, p)(\overline{p}-x_0)u,(K_\Phi(\cdot, p)-K_\Phi(\cdot,
x_0))u\right).
\label{reatv1}
\end{equation}
Then the following holds:
\begin{equation}
\label{qwertyu}
(f,g)\in \mathcal R_{x_0}\quad\Longrightarrow\quad [f,f]=[f+2x_0g,f+2x_0g].
\end{equation}

We first prove that
\begin{equation}
\label{qwertyu1} [f,g]+[g,f]+2x_0[g,g]=0.
\end{equation}

An element in $\mathcal R_{x_0}$ can be written as $(f,g)$ with
\begin{equation}
\label{f-and-g}
\begin{split}
f(p)&=\sum_{j=1}^mK_\Phi(p,p_j)(\overline{p_j}-x_0)u_j\\
g(p)&=\sum_{j=1}^mK_\Phi(p,p_j)u_j-K_\Phi(p,x_0)d,\quad {\rm
where}\quad d=\sum_{j=1}^m u_j.
\end{split}
\end{equation}
With $f$ and $g$ as in \eqref{f-and-g} we have:
\[
\begin{split}
[f,g]&=\left(\sum_{i,j=1}^m
u_i^*K_\Phi(p_i,p_j)(\overline{p_j}-x_0)u_j\right)-
d^*\left(\sum_{j=1}^mK_\Phi(x_0,p_j)(\overline{p_j}-x_0)u_j\right),\\
[g,f]&=\left(\sum_{i,j=1}^m
u_i^*(p_i-x_0)K_\Phi(p_i,p_j)u_j\right)-
\left(\sum_{i=1}^mu_i^*(p_i-x_0)K_\Phi(p_i,x_0)\right)d.
\end{split}
\]
Thus
\[
\begin{split}
[f,g]+[g,f]&=-2x_0\left(\sum_{i,j=1}^m
u_i^*K_\Phi(p_i,p_j)u_j\right)+\\
&\hspace{5mm}+\sum_{i,j=1}^mu_i^*\left\{p_iK_\Phi(p_i,p_j)+K_\Phi(p_i,p_j)
\overline{p_j}\right\}u_j-d^*\left(\sum_{j=1}^mK_\Phi(x_0,p_j)\overline{p_j}u_j
\right)+\\
&\hspace{5mm}+x_0d^*\left(\sum_{j=1}^mK_\Phi(x_0,p_j)u_j\right)-
\left(\sum_{i=1}^mu_i^*p_iK_\Phi(p_i,x_0)\right)d+\\
&\hspace{5mm}+x_0\left(\sum_{j=1}^mu_j^*K_\Phi(p_j,x_0)\right)d.
\end{split}
\]
Taking into account \eqref{equaphi} we have
\[
\begin{split}
[f,g]+[g,f]&=-2x_0\left(\sum_{i,j=1}^m
u_i^*K_\Phi(p_i,p_j)u_j\right)+\\
&\hspace{5mm}+\left(\sum_{i=1}^mu_i^*J\Phi(p_i)\right)d+d^*
\left(\sum_{j=1}^m\Phi(p_j)^*Ju_j\right)-\\
&\hspace{5mm}-d^*\left(\sum_{j=1}^mK_\Phi(x_0,p_j)\overline{p_j}u_j
\right)+x_0d^*\left(\sum_{j=1}^mK_\Phi(x_0,p_j)u_j\right)-\\
&\hspace{5mm}-\left(\sum_{i=1}^mu_i^*p_iK_\Phi(p_i,x_0)\right)d+x_0
\left(\sum_{i=1}^mu_i^*K_\Phi(p_i,x_0)\right)d.
\end{split}
\]
We now turn to $[g,g]$. We have:
\[
\begin{split}
[g,g]&=\left(\sum_{i,j=1}^mu_i^*(K_\Phi(p_i,p_j)u_j\right)-d^*
\left(\sum_{j=1}^mK_\Phi(x_0,p_j)u_j\right)-\\
&\hspace{5mm}-\left(\sum_{i=1}^mu_i^*K_\Phi(p_i,x_0)\right)d+d^*K_\Phi(x_0,x_0)d.
\end{split}
\]
Thus
\[
\begin{split}
[f,g]+[g,f]+2x_0[g,g]&=
\left(\sum_{i=1}^mu_i^*J\Phi(p_i)\right)d+d^*
\left(\sum_{j=1}^m\Phi(p_j)^*Ju_j\right)-\\
&\hspace{5mm}-d^*\left(\sum_{j=1}^mK_\Phi(x_0,p_j)\overline{p_j}u_j
\right)-x_0d^*\left(\sum_{j=1}^mK_\Phi(x_0,p_j)u_j\right)-\\
&\hspace{5mm}-\left(\sum_{i=1}^mu_i^*p_iK_\Phi(p_i,x_0)\right)d-
x_0\left(\sum_{i=1}^mu_i^*K_\Phi(p_i,x_0)\right)d+\\
&\hspace{5mm}+2x_0d^*K_\Phi(x_0,x_0)d\\
&= \left(\sum_{i}^mu_i^*J\Phi(p_i)\right)d+d^*
\left(\sum_{j=1}^m\Phi(p_j)^*Ju_j\right)-\\
&\hspace{5mm}-d^*\left(\sum_{j=1}^mK_\Phi(x_0,p_j)(\overline{p_j}+x_0)u_j
\right)-\\
&\hspace{5mm}-\left(\sum_{j=1}^mu_i^*(p_i+x_0)K_\Phi(p_i,x_0)\right)d+
2d^*x_0K_\Phi(x_0,x_0)d
\end{split}
\]
using
\[
K_\Phi(x_0,x_0)=\frac{1}{2x_0}\left(J\Phi(x_0)+\Phi(x_0)^*J\right),
\]
we obtain
\[
\begin{split}
[f,g]+[g,f]+2x_0[g,g]&=
\left(\sum_{i}^mu_i^*J\Phi(p_i)\right)d+d^*
\left(\sum_{j=1}^m\Phi(p_j)^*Ju_j\right)-\\
&\hspace{5mm}-d^*\left(\sum_{j=1}
   \left(J\Phi(x_0)+\Phi(p_j)^*J\right)u_j\right)-\\
&\hspace{5mm}-\left(\sum_{j=1}^mu_i^*
    \left(J\Phi(p_i)+\Phi(x_0)^*J\right)\right)d+
2d^*x_0K_\Phi(x_0,x_0)d\\
&=0
\end{split}
\]
and so we have proved \eqref{qwertyu1}. Equation \eqref{qwertyu}
follows since
\[
[f+2x_0g,f+2x_0g]=[f,f]+2x_0\left([f,g]+[g,f]+2x_0[g,g]\right).
\]

Equation \eqref{qwertyu1} expresses that the linear space of
functions $(f,f+2x_0g)$ with $f,g$ as in \eqref{f-and-g} define an
isometric relation $\mathcal R$ from the Pontryagin space
$\mathcal L(\Phi)$ into itself. Let now $h\in\mathcal L(\Phi)$ be
such that
\[
[h,K_\Phi(\cdot,p)(\overline{p}-x_0)u]= 0 \ \ \forall
p\in\Omega\quad{\rm and}\quad u\in\mathcal H.
\]
Then
\[
(p-x_0)h(p)=0,\quad \forall p\in\Omega
\]
and $h\equiv 0$ in $\Omega$ (recall that the elements of
$\mathcal L(\Phi)$ are slice hyperholomorphic in $\Omega$). Thus
the domain of this relation is dense. By the quaternionic version
of Shmulyan's theorem (see \cite[Theorem 7.2]{acs2}), $\mathcal R$
extends to the graph of a (continuous) isometry, say $B^*$, on
$\mathcal L(\Phi)$. We have for $h\in\mathcal L(\Phi)$
\[
\begin{split}
u^*(p-x_0)((Bh(p))&=[Bh,K_\Phi(\cdot, p)(\overline{p}-x_0)u]\\
&=[h, B^*(K_\Phi(\cdot, p)(\overline{p}-x_0)u)]\\
&=[h,K_\Phi(\cdot, p)(\overline{p}-x_0)u+2(K_\Phi(\cdot, p)-K_\Phi(\cdot, p))u]\\
&=u^*\left((p-x_0)h(p)+2h(x_0)-2h(x_0)\right)\\
&=u^*\left((p+x_0)h(p)-2h(x_0)\right).
\end{split}
\]

We note that $\mathcal R_{x_0}$ extends to the graph of $R_{x_0}^*$.\\

STEP 2: {\sl The function $p\mapsto R_{x_0}\Phi \eta$ belongs to
$\mathcal L(\Phi)$ for every $\eta\in\mathcal H$ and the operator $F$
from $\mathcal H$ into $\mathcal L(\Phi)$ defined by
\[
F\eta=R_{x_0}\Phi\eta
\]
is bounded.}\\

We note that $B=I+2x_0R_{x_0}$ and so $R_{x_0}$
is a bounded operator in $\mathcal L(\Phi)$. From \eqref{equaphi}
we have for $\xi\in\mathcal H$
\begin{equation}
J\Phi(p)\xi+\Phi(x_0)^*J\xi=pK_\Phi(p,x_0)\xi+K_\Phi(p,x_0)\xi
x_0.
 \label{asdfg}
\end{equation}
Apply $R_{x_0}$ on both sides (as an operator on slice
hyperholomorphic functions; the two sides of \eqref{asdfg} will
no belong to $\mathcal L(\Phi)$ in general). Note that
\[
R_{x_0}(pf(p))=f(p)+x_0(R_{x_0}f)(p),
\]
and so we obtain
\[
R_{x_0}\Phi
J\xi=K_\Phi(p,x_0)\xi+x_0(R_{x_0}K_\Phi(\cdot,x_0)\xi)(p)+(R_{x_0}K_\Phi(p,x_0)\xi
x_0)(p),
\]
and this expresses $R_{x_0}\Phi J\xi$ as an element of $\mathcal
L(\Phi)$ since $R_{x_0}$ is bounded in $\mathcal L(\Phi)$ and so
the elements on the right side of the above equality belong to
$\mathcal L(\Phi)$. This ends the proof of the first claim since $J$ is invertible. Finally, to see that the operator
$F$ is bounded we remark that it is closed and everywhere defined.\\

We remark that the argument is different from the one for the corresponding operator $F$
(defined by \eqref{I2}) in the spaces $\mathcal P(S)$. In the classical case, the argument we are aware of, uses
a Cayley transform to go back to the case of generalized Schur functions. The argument we presented here is
probably known in the classical case, but we are not aware of any reference for it.\\

STEP 3: {\sl The realization formula \eqref{eq:realphi12345} holds.}\\

The proof is the same as the one in STEP 6 for $S$.\\

STEP 4: {\sl The function $\Phi$ admits a slice hypermeromorphic
extension to $\mathbb H_+$}\\

Recall that $T=I+2x_0B$ is co-isometric. Using Theorem \ref{tm:inv} we can thus write $T$ as
\[
T=\begin{pmatrix}T_{11}&T_{12}\\0&T_{22}\end{pmatrix}
\]
where $T_{11}$ is a bijective contraction from a anti-Hilbert
space onto itself, and $T_{22}$ is a contraction from a Hilbert
space into itself. Thus for $x>0$ in a neighborhood of $x_0$,
\[
\begin{split}
(1+x)I+(x-1)B&=(1+x)I+(x-1)\left(\frac{T-I}{2}\right)\\
&=(3+x)+(x-1)T\\
&=(3+x)\left(I+\begin{pmatrix}\frac{x-1}{3+x}T_{11}&\frac{x-1}{3+x}T_{12}\\0&
\frac{x-1}{3+x}T_{22}\end{pmatrix}\right)
\end{split}
\]
and hence the result by slice hyperholomorphic extension since $\frac{q-1}{3+q}$ sends $\mathbb H_+$
into $\mathbb B_1$.\\

STEP 5: {\sl A function $\Phi$ admitting a realization of the form
\eqref{eq:realphi12345} is in a class ${\rm
GP}_\kappa(J)$.}\smallskip

The proof is as in the case of the functions $S$ and is based on
the identity
\[
J\Phi(x)+\Phi(y)^*J=(x+y)G(I(x_0+x)-(x_0-x)B)^{-1}(I(x_0+y)-(x_0-y)B)^{-*},
\]
where $x,y$ are real and in a neighborhood of $x_0$.\\

STEP 6: {\sl An observable realization of the form
\eqref{eq:realphi12345} is unique up to a isomorphism of quaternionic Pontryagin
spaces.}\\
\end{proof}

We note that the relation (\ref{reatv1}) is inspired from \cite[p.
708]{atv1} and more generally, by the constructions of the
"$\epsilon$-method" developed in the papers of Krein and Langer;
see for instance \cite{kl1,MR47:7504} for the latter.

\begin{Cy}
In $\mathcal L(\Phi)$ it holds that
\begin{equation}
\label{r1r1}
R_{x_0}+R_{x_0}^*=-2x_0R_{x_0}^*R_{x_0}.
\end{equation}
\end{Cy}

\begin{proof}
This is a rewriting of \eqref{qwertyu}.
\end{proof}

We note that \eqref{r1r1} is a special case of the structural
identity characterizing $\mathcal L(\Phi)$ spaces in the complex
case, see \cite{dbbook}. To ease the notation we consider the
case $J=I$.

\begin{Tm}
Let the $\mathbf L(\mathcal H)$-valued function $\Phi$ be slice
hyperholomorphic in an axially symmetric s-domain  $\Omega$
containing $p=0$, such that the associated space does not contain
non zero constants, and has its elements slice hyperholomorphic
in a neighborhood of the origin. Assume that $\Phi\in{\rm
GP}_\kappa(I)$. Then there exists a right quaternionic Hilbert
space $\mathcal H_1$ and operators
\[
\begin{pmatrix}A&B\\ C&D\end{pmatrix}:\,\,\,
{\mathcal H_1}\oplus\mathcal H\,\,\longrightarrow\,\, {\mathcal
H_1}\oplus\mathcal H
\]
such that
\begin{equation}
\label{realphi!!!}
\Phi(p)=D+p C\star (I_{\mathcal
H_1}-pA)^{-\star} B
\end{equation}
and
\[
{\rm Re}\,\begin{pmatrix}A&B\\
C&D\end{pmatrix}\begin{pmatrix}I_{\mathcal H_1}&0\\0&-I_{\mathcal
H}\end{pmatrix}=0.
\]
\end{Tm}

\begin{proof}

We first define a linear relation $R_\Phi$ in $(\mathcal
L(\Phi)\oplus \mathcal H)\times (\mathcal L(\Phi)\oplus \mathcal
H)$ via the formulas

\begin{equation}
\label{realJ}
\left(\begin{pmatrix} -K_\Phi(\cdot, q)\overline{q}u\\
                      u\end{pmatrix}\,,\,
                      \begin{pmatrix} K_\Phi(\cdot, q)u
                      \\
                     \Phi(q)^* u
                      \end{pmatrix}
                      \right).
\end{equation}

STEP 1: {\sl The relation $R_\Phi$ satisfies
\begin{equation}
\label{eqr2} {\rm Re}\,
\langle  \begin{pmatrix}f\\ -g\end{pmatrix}\,,\,\begin{pmatrix}F\\
G\end{pmatrix}\rangle= 0.
\end{equation}
Furthermore, it has dense domain since the space $\mathcal
L(\Phi)$
contains no non zero constant functions.}\\
Let
\[
\begin{pmatrix}f\\ -g\end{pmatrix}=-\sum_{n=1}^t \begin{pmatrix} K_\Phi(\cdot,
q_n)\overline{q}_n
u_n\\
                      u_n\end{pmatrix}\,,\,
\]
and
\[
\begin{pmatrix}F\\ G\end{pmatrix}=\sum_{n=1}^t
\begin{pmatrix} K_\Phi(\cdot, q_n)u_n
                      \\
                     \Phi(q_n)^* u_n
                      \end{pmatrix}.
\]
Then
\[
\begin{split}
\langle\begin{pmatrix}f\\ -g\end{pmatrix},\begin{pmatrix}F\\
G\end{pmatrix}\rangle&=-\sum_{n,m=1}^t
u_m^*K_\Phi(q_m,q_n)\overline{q_n}u_n+u_m^*\Phi(q_n)u_m
\end{split}
\]
so that, using \eqref{equaphi}, we obtain
\[
\begin{split}
{\rm Re}\,\langle\begin{pmatrix}f\\ -g\end{pmatrix},\begin{pmatrix}F\\
G\end{pmatrix}\rangle&=0.
\end{split}
\]

STEP 2: {\sl The relation $R_\Phi$ is the graph of a densely
defined operator which has a continuous extension, and its
adjoint is the backward-shift realization
\[
\begin{pmatrix}A&B\\ C&D\end{pmatrix},
\]
where
\[
\begin{split}
pAf(p)&=f(p)-f(0),\\
pBu(p)&=(\Phi(p)-\Phi(0))u,\\
Cf&=f(0),\\
Du&=\Phi(0)u.
\end{split}
\]
}

We only have to consider the operator $B$. Consider a family $T$
of pairs $(q,u)\in\Omega\times \mathcal H$ such that the
functions $K_\Phi(\cdot, q)u$ are linearly independent and span
the space of all the functions $K_\Phi(\cdot, p)v$, where $p$
runs through all of $\Omega$ and $v$ runs through all of
$\mathcal H$. Define a densely defined operator from $\mathcal
L(\Phi)$ into $\mathcal H$ by
\[
X(K_\Phi(\cdot, p)u)=(\Phi(p)^*-\Phi(0)^*)u,\quad (p,u)\in A.
\]

We claim that $X$ has an adjoint which is the operator $B$ above.
To see that, we remark that \eqref{eqr2} can be rewritten as

\begin{equation}
\label{eqr3} \langle  \begin{pmatrix}f\\
-g\end{pmatrix}\,,\,\begin{pmatrix}
A^*&C^*\\ X&D^*\end{pmatrix}\begin{pmatrix}f\\
g\end{pmatrix}\rangle+ \langle \begin{pmatrix}
A^*&C^*\\ X&D^*\end{pmatrix}\begin{pmatrix}f\\
g\end{pmatrix}\,,\, \begin{pmatrix}f\\
-g\end{pmatrix}\rangle\ =0.
\end{equation}

Using the quaternionic polarization formula it follows that for
any
\[
\begin{pmatrix}f_1\\
g_1\end{pmatrix}\quad{\rm and}\quad \begin{pmatrix}f_2\\
g_2\end{pmatrix}
\]
in the domain of $R_\Phi$ we have

\begin{equation}
\label{eqr34} \langle  \begin{pmatrix}f_1\\
-g_1\end{pmatrix}\,,\,\begin{pmatrix}
A^*&B*\\ X&D^*\end{pmatrix}\begin{pmatrix}f_2\\
g_2\end{pmatrix}\rangle+ \langle \begin{pmatrix}
A^*&B^*\\ X&D^*\end{pmatrix}\begin{pmatrix}f_1\\
g_1\end{pmatrix}\,,\, \begin{pmatrix}f_2\\
-g_2\end{pmatrix}\rangle\,\,=0
\end{equation}
and so $R_\Phi$ has an adjoint and so does $X$.

\end{proof}
It is useful to note that the operator $B$ appearing in the previous theorem is the opposite of the operator in (\ref{I1}).
\begin{Ex}{\rm
As an illustration of the previous theorem consider the function
\begin{equation}
\label{asdf} \varphi(p)=(p+a)^{-\star},
\end{equation}
where $a\in\mathbb H$ is such that ${\rm Re}\, a=0$. Set
\[
M(p,q)=(p+a)^{-\star}\overline{(q+a)^{-\star}}.
\]
Since $a+\overline{a}=0$ we have
\[
p
M(p,q)+M(p,q)\overline{q}=\varphi(p)+\overline{\varphi(q)},
\]
and so $\varphi$ is a positive function. For $p=x>0$ we have
\[
\varphi(x)=a^{-1}-\frac{x}{(1+ xa^{-1})a^2},
\]
which leads to the realization \eqref{realphi!!!} with
\[
\begin{pmatrix}A&B\\ C&D\end{pmatrix}=\begin{pmatrix}-a^{-1}&a^{-1}\\
-a^{-1}&a^{-1}\end{pmatrix}
\]
So
\[
{\rm Re}\,\begin{pmatrix}A&B\\
C&D\end{pmatrix}\begin{pmatrix}1&0\\0&-1\end{pmatrix}
=-(a^{-1}+\overline{a^{-1}})
\begin{pmatrix}
1&1\\1&1\end{pmatrix}=\begin{pmatrix}0&0\\0&0\end{pmatrix}
\]
since $a+\overline{a}=0$.\\
}
\end{Ex}

\subsection{The positive case}

In this section we prove one theorem in the case $\kappa=0$ and
$J=I$. We say that the function $\Phi$ is positive rather that
writing $\Phi\in{\rm GP}_0(I)$. The proof uses the existence of a
squareroot of a positive operator in a quaternionic Pontryagin
space. In the indefinite case, such a result still exists in the
complex case (this is called the Bognar-Kramli theorem, see
\cite[Theorem 2.1 p. 149]{bognar}, \cite[Theorem
1.1.2]{MR92m:47068}). A quaternionic version of this
factorization theorem is not available at present.

\begin{Tm}
Let $\Phi$ be slice-hyperholomorphic in an axially symmetric
s-domain of the origin with realization \eqref{realphi!!!} such
that
\[
{\rm Re}\,\begin{pmatrix}A&B\\
C&D\end{pmatrix}\begin{pmatrix}I_{\mathcal H_1}&0\\0&-I_{\mathcal
H}\end{pmatrix}\le 0
\]
Then $\Phi$ is positive.
\end{Tm}
\begin{proof}

We first note that a positive operator $T$, in a quaternionic
Hilbert space has a squareroot, that is, there exists a positive
operator $X$ such that $X^2=T$. The proof uses the spectral
theorem, which holds for Hermitian operators in quaternionic
Hilbert spaces. The theorem is mentioned without proof in a
number of papers (see for instance \cite{MR0137500},
\cite{viswanath_thesis}, \cite{MR44:2067}). The spectrum used in
these works is not the $S$-spectrum, see \cite[p.
141]{MR2752913}); a proof is given in the preprint \cite{acs4}.
Another way to prove the existence of a squareroot is to define
(assuming first $\|T\|\le 1$), as in the complex case, a sequence
of operators $X_0,X_1,\ldots$ by $X_0=0$ and
\[
X_{n+1}=\frac{1}{2}((I-T)+X_n^2),\quad n=0,1,\ldots,
\]
(see for instance \cite[p. 64]{MR675952}) and check that:\\
$(1)$ A weakly convergent increasing sequence of
positive operators converges strongly.\\
$(2)$ An increasing family $(X_n)_{n\in\mathbb N}$ of bounded
positive operators such that
\[
\lim_{n\rightarrow\infty} \langle
X_nf,f\rangle<\infty,\quad\forall f\in\mathcal H
\]
converges strongly to a positive operator. Since the arguments do
not differ from the complex case we omit them.\\

Let $X$ be the squareroot of $-{\rm Re}\,\begin{pmatrix}A&B\\
C&D\end{pmatrix}\begin{pmatrix}I_{\mathcal H_1}&0\\0&-I_{\mathcal
H}\end{pmatrix}$. We write
\[
X=\begin{pmatrix} L\\ K\end{pmatrix}
\]
where $L$ is a linear operator from $\mathcal H\times \mathcal
H_1$ into $\mathcal H$ and $K$ is a a linear operator from
$\mathcal H\times \mathcal H_1$ into $\mathcal H_1$.\\

Let now
\[
{\rm Re}\,\begin{pmatrix}A&B\\
C&D\end{pmatrix}\begin{pmatrix}I_{\mathcal H_1}&0\\0&-I_{\mathcal
H}\end{pmatrix}= -\begin{pmatrix} L\\ K\end{pmatrix}
\begin{pmatrix}
L\\ K\end{pmatrix}^*,
\]
Then

\[
\begin{split}
\Phi(x)+\Phi(y)^*&=D+D^*+xC(I-xA)^{-1}B+yB^*(I-yA)^{-*}C^*\\
&=KK^*+xC(I-xA)^{-1}(C^*-LK^*)+\\
&\hspace{5mm}+y(C-KL^*)(I-yA)^{-*}C^*\\
&=(K-xC(I-xA)^{-1}L)(K-yC(I-yA)^{-1}L)^*+\\
&\hspace{5mm}+xC(I-xA)^{-1}C^*+yC(I-yA)^{-*}C^*-\\
&\hspace{5mm}-xyC(I-xA)^{-1}LL^*(I-yA)^{-*}C^*.
\end{split}
\]
But, using $A+A^*+LL^*=0$, we have
\[
\begin{split}
xC(I-xA)^{-1}C^*+yC(I-yA)^{-*}C^*-
xyC(I-xA)^{-1}LL^*(I-yA)^{-*}C^*&=\\
(x+y)C(I-xA)^{-1}(I-yA)^{-*}C^*
\end{split}
\]
The claim follows by slice-hyperholomorphic extension.
\end{proof}

We note that the computations are classical, see for instance
\cite{MR525380}, \cite[Theorem 3.3, p. 26]{faurre}.\\

We conclude this section with an example of elements of ${\rm GP}_0(I)$ (that is, positive functions) which play an
important role in models for pairs of anti self-adjoint operators.
This originates with the paper of de Branges and Rovnyak
\cite{dbr1}. We refer to \cite{MR1960423,MR2069002,Iacob} for
examples and applications of the model of de Branges and Rovnyak.
In this section, we briefly outline how a positive function also
appears in the present setting. We follow the approach of
\cite{MR2069002}, and consider bounded operators for the sake of
illustration. The proof of the following lemma is as in \cite[p.
18]{MR2069002} and is omitted.

\begin{Pn}
\label{pndbr}
Let $T_+$ and $T_-$ be two anti-self-adjoint operators in the
quaternionic space $\mathcal H$. Then:\\
$(1)$ The space
\begin{equation}
\label{inv}
\cap_{u=1}^\infty \ker (T_+^u-T_-^u)
\end{equation}
is the largest subspace, invariant under $T_+$ and $T_-$ and on
which they coincide.\\
$(2)$ Assume that ${\rm rank}\, T_+-T_-=n<\infty$. Then there
exists a $n\times n$ matrix $J\in\mathbb H^{n\times n}$ such that
$J^2=-I_n$ and $J^*=-J$, and a linear bounded operator $C$ from
$\mathcal H$ into $\mathbb H^n$ such that
\[
T_+-T_-=-C^*JC.
\]
\end{Pn}

\begin{Tm} Using the notation of the preceding lemma,
the function
\[
\Phi(p)=J+C\star (pI-T_+)^{-\star}C^*
\]
is positive and its inverse is equal to
\[
\Phi^{-\star}(p)=-J-JC\star(pT-T_-)^{-\star}C^*J.
\]
\end{Tm}

\begin{proof}
For $x,y$ on the positive real axis we have (recall that $T_+^*=-T_+$)
\[
\begin{split}
\Phi(x)+\Phi(y)^*&=J+J^*+C(xI-T_+)^{-1}C^*+C(yI-T_+)^{-*}C^*)\\
&=C(xI-T_+)^{-1}C^*+C(yI+T_+)^{-1}C^*\\
&=xK(x,y)+K(x,y)y,
\end{split}
\]
where $K(x,y)=C(xI-T_+)^{-1}(yI-T_+)^{-*}C^*$. The result follows then by slice hyperholomorphic extension.\smallskip

Still for positive $x$ and using for instance formula \eqref{formreal1} we have
\[
\begin{split}
\Phi(x)^{-1}&=J^{-1}-J^{-1}C(xI-(T_+-C^*J^{-1}C^*)^{-1}C^*J^{-1}\\
&=-J-JC(xI-T_-)^{-1}C^*J.
\end{split}
\]
The formula for $\Phi^{-1}$ follows then by slice hyperholomorphic extension.
\end{proof}

When the space \eqref{pndbr} in Proposition \ref{pndbr} is trivial the function $\Phi$ characterizes the pair $(T_+,T_-)$.
Models for pairs of (possibly unbounded) anti-self-adjoint
operators in a quaternionic Hilbert space in terms of the
reproducing kernel Hilbert spaces $\mathcal L(\Phi)$ and
$\mathcal L(\Phi^{-1})$, and related trace formulas similar to
the ones presented in the papers \cite{dbr1,MR2069002} will be
considered elsewhere.

{\bf Acknowledgments}. The authors are grateful to the anonymous referee for carefully reading the manuscript and for the useful comments.

\bibliographystyle{plain}
\def\cprime{$'$} \def\lfhook#1{\setbox0=\hbox{#1}{\ooalign{\hidewidth
  \lower1.5ex\hbox{'}\hidewidth\crcr\unhbox0}}} \def\cprime{$'$}
  \def\cfgrv#1{\ifmmode\setbox7\hbox{$\accent"5E#1$}\else
  \setbox7\hbox{\accent"5E#1}\penalty 10000\relax\fi\raise 1\ht7
  \hbox{\lower1.05ex\hbox to 1\wd7{\hss\accent"12\hss}}\penalty 10000
  \hskip-1\wd7\penalty 10000\box7} \def\cprime{$'$} \def\cprime{$'$}
  \def\cprime{$'$} \def\cprime{$'$}

\end{document}